\documentclass{article}

\usepackage{amsmath}
\usepackage{amssymb}
\usepackage{graphicx}

\begin{document}

\renewcommand{\theequation}{\arabic{section}.\arabic{equation}}

\bigskip

\begin{center}
{\Large Identification and control of SARS-CoV-2 epidemic model parameters }

\medskip

{\Large Gabriela Marinoschi}

\medskip

\textquotedblleft Gheorghe Mihoc-Caius Iacob\textquotedblright\ Institute of
Mathematical Statistics and

Applied Mathematics of the Romanian Academy,

Calea 13 Septembrie 13, Bucharest, Romania

gabriela.marinoschi@acad.ro

\medskip 
\end{center}

\noindent{\small Abstract. We propose a mathematical model with five
compartments for the SARS-CoV-2 transmission: susceptible }${\small S},$%
{\small \ undetected infected asymptomatic }${\small A},${\small \
undetected infected symptomatic }${\small I}${\small , confirmed positive
and isolated }${\small L},${\small \ and recovered }${\small R},${\small \
for which we have a twofold objective. First ,we formulate and solve an
inverse problem focusing mainly on the identification of the values }$%
{\small A}_{0}${\small \ and }${\small I}_{0}${\small \ of the undetected
asymptomatic and symptomatic individuals, at a time }$t_{0},${\small \ by
available measurements of the isolated and recovered individuals at two
succeeding times, }${\small t}_{0}${\small \ and }${\small T>t}_{0}.${\small %
\ Simultaneously, we identify the rate standing for the average number of
individuals infected in unit time by an infective symptomatic individual.
Then, we propose a control problem aiming at controlling the infected
classes by improving the actions in view of isolating as much as possible
the populations }${\small A}${\small \ and }${\small I}${\small \ in the
class }${\small L.}${\small \ These objectives are formulated as
minimization problems, the second one including a state constraint, which
are treated by an optimal control technique. The existence of optimal
controllers is proved and the first order necessary conditions of optimality
are determined. For the second problem, they are deduced by passing to the
limit in the conditions of optimality calculated for an appropriately
defined approximating problem. In this case, the dual system is singular and
has a component in the space of measures. The discussion of the asymptotic
stability of the system done for the case when life immunity is gained
reveals an asymptotic extinction of the disease, with a well determined
reproduction number.}

\medskip


\medskip


\section{Introduction}

\setcounter{equation}{0}

The current pandemics of COVID-19 disease caused by the severe acute
respiratory syndrome coronavirus 2 (SARS-CoV-2) has undergone an accentuated
exponential increase of cases all over the world. As in other transmissible
diseases, the infected with SARS-CoV-2 may have in the incubation period
mild forms or even no symptoms such that they can be not aware of the fact
that are carrying the virus. They are known as exposed. But, the particular
and the worst aspect of COVID-19 disease is that the exposed, called here
asymptomatic, are highly contagious and can transmit the disease (see \cite%
{Cascela}). Hence, the early identification of individuals infected with
SARS-CoV-2 and the necessity of isolation of the people found infected is
crucial for reducing the virus spread. Mathematical modeling can help to
estimate some relevant parameters of this epidemic, which can allow the
prediction of the disease evolution and the preparation of the necessary
measures for the disease containment.

The study of various aspects of the SARS-CoV-2 has led to an extremely rich
article production since the debut of the pandemics. Mathematical modelling
of various aspects of the disease has been addressed as well. Many of them
are based on the SIR model which describes the transmission of the disease
through three stages of infection, susceptible, infected and recovered.
However, taking into account the previous considerations, the SIR model
cannot adequately characterize the COVID-19 disease. Several more complex
models have been also considered to illustrate the spread of this disease.
We cite here a few studies: a SEIR (susceptible, exposed, infectious,
removed) model considering risk perception and the cumulative number of
cases has been developed in \cite{Lin}; a discrete-time SIR model including
dead individuals was proposed in \cite{Anast}, a control-oriented SIR model
that puts into evidence the effects of delays and compares the outcomes of
different containment policies was discussed in \cite{Casella}. In \cite%
{Pu-So} a mathematical method was developed to deduce the evolution over
time of the new coronavirus infection and to establish the effect of
isolation strategies from the accumulated data, such as the number of deaths
and hospitalizations. A more detailed model of transmission in Italy that
extends the classic SEIR model was presented and analyzed in \cite{Gio}.
This model, called SIDARTHE, involves many compartments, such as:
susceptible, non-life-threatening cases, asymptomatic with minor and
moderate infection, symptomatic, for each of them being separate classes of
detected and undetected individuals, symptomatic with a severe situation,
dead and recovered. The model omits the probability rate of becoming
susceptible again after having recovered from the infection. The model
parameters were estimated by a best-fit approach, namely by finding the
parameters that locally minimize the sum of the squares of the errors. The
computations were based on data measured in Italy between the beginning of
the outbreak and early April and were updated over time to reflect the
progressive introduction of increased restrictions.

In this paper we introduce a mathematical model for SARS-CoV-2 epidemic,
involving five compartments considered to be essential to depict the feature
of the epidemic. A first goal is to approach by an optimal control technique
the identification of some parameters such as the average number of
individuals infected in unit time by an infected symptomatic and the number
of asymptomatic and symptomatic individuals still undetected at a moment of
time. These parameters can be further used to calculate the reproduction
rate and are relevant for predicting the evolution of the epidemic and for
planning effective control policies. These justifies the second objective in
the paper which refers to the control of classes $A$ and $I,$ by finding
optimal control coefficients related to actions as screening, testing,
tracing, which may lead to the reduction of the infected population by
moving it into the isolated compartment $L$.

\subsection{Mathematical model}

We introduce a mathematical model of SEIR type with five compartments,
represented at time $t$ by: susceptible $S(t)$ (a healthy individual which
can acquire the disease), infected asymptomatic $A(t)$ (individuals who have
acquired the disease, have no symptoms or only mild ones and have been not
tested yet)$,$ infected with symptoms but not confirmed yet $I(t),$ detected
by tests and isolated (by hospitalization, quarantine, isolation) $L(t)$ and
recovered $R(t).$ As said before, the individuals in the infected classes $A$
and $I,$ meaning those not tested yet, and whose number is not known, can be
in circulation and infect other people. The people tested and found
infectious are supposed to be removed from circulation until healing and
introduced in the class $L$. Since the identification of the parameters we
propose here is considered to be performed on short intervals of time, we
skip in this model the extinct (from any reason) and the newborns, because
their number is negligible compared with the large number of people in the
relevant compartments for this disease.

$$\begin{tabular}{c}  
\includegraphics[scale=.7]{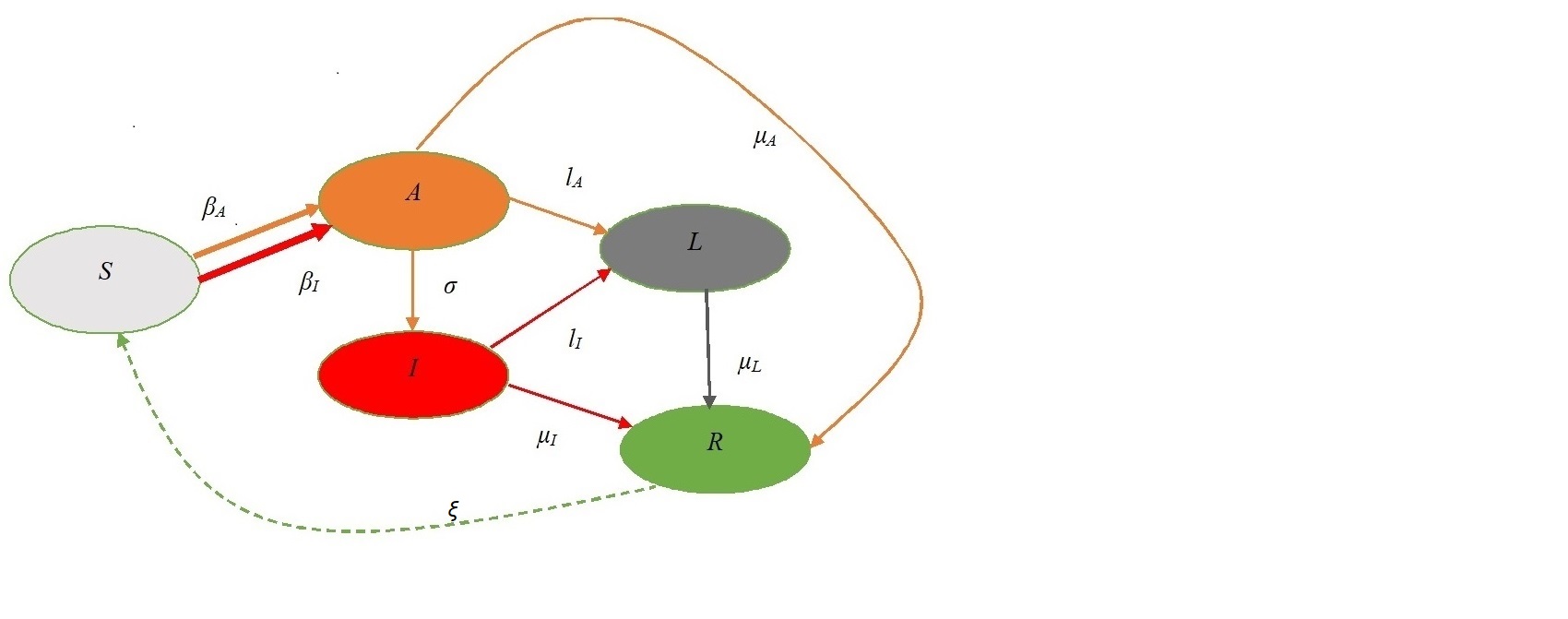}\label{Fig1}\\
\hspace*{-80mm}
{\bf Figure 1.  Flowchart of the SAILR model}
\end{tabular}$$

Thus, we assume that both $A$ and $I$ are infectious and so the susceptible
can be infected by the infective with symptoms with the rate $\beta _{I}(t),$
as well as by the infective asymptomatic with a rate $\beta _{A}(t).$ We
also assume that the infective $A$ and $I$ detected by various actions are
quickly isolated, in order to help the disease containment. Then, the
individuals $A$ have three possibilities: can become symptomatic with the
rate $\sigma ,$ can be detected by testing and isolated with the rate $l_{A}$
or can recover with the rate $\mu _{A}.$ At their turn, the infected
symptomatic can be detected by testing and isolated with the rate $l_{I}$
and can recover with the rate $\mu _{I}.$ The individuals isolated recover
with the rate $\mu _{L}.$ Existing evidence (see e.g., \cite{Linton}) shows
that the immunity is gained for a short time only, so that the recovered can
lose it after some time and go back into the susceptible compartment with
the rate $\xi (t).$ Of course, the model can be completed with other classes
of population but in this paper we keep only this formulation. We call the
model SAILR and depict its flowchart in Fig. 1.

The model is going be written for fractions of the total population (that is 
$S,A,I,L,R$ represent the real number of the individuals in these classes
divided by the total population $N(t)).$ For the SIR model the form of the
normalized system is described in the monograph \cite{Ia-Pu}. There, $\beta
_{I}$ stands for\ the average number of individuals infected in unit time by
an infective symptomatic. This is given by the number of contacts an
infective $I$ has in the time unit, multiplied by the probability that a
contact produces an infective, when one of the two individuals is
susceptible and the other is infective symptomatic. For the SAILR model, we
derive in a similar way as in \cite{Ia-Pu} the normalized equations, by
simple calculations, but we do no longer write them. Thus, we can
analogously characterize $\beta _{A}$ as being the average number of
individuals infected in unit time by an infective asymptomatic. The rates $%
\beta _{I},\beta _{A}$ and $\xi $ can vary in time. For example, $\xi $ can
be zero for some time and then begin to increase, and $\beta _{I}$ and $%
\beta _{A}$ may have a periodic increasing-decreasing behavior due to the
particularities of the transmission, the medical or social measures that are
imposed or the variation of the virus virulence. The removal rates $\mu _{I},
$ $\mu _{A}$ and also $\mu _{L}$ are considered constant because they
generally depend on the interaction between the pathogen agents and the
immune system of an infected individual (see e.g. \cite{Ia-Pu}). The
parameters $\tau _{I}=1/\mu _{I},$ $\tau _{A}=1/\mu _{A}$ and $\tau
_{L}=1/\mu _{L}$ are the average durations of the infection for an infective 
$I,$ an infective $A$ (which may follow a treatment or not) and an infective
isolated $L$, respectively. The last one is supposed to receive a treatment.

The rates $l_{A}$ and $l_{I}$ are related to the probability rate of
detection, relative to asymptomatic and symptomatic cases, respectively.
They may reflect for instance the number of tests performed over the
population and they can be modified by enforcing sustained actions, as for
example a massive testing campaign (see \cite{Peto}). The value $l_{I}$ may
be larger than $l_{A}$, because a symptomatic individual is more likely to
be tested. We assume that $l_{I}$ and $l_{A}$ are constant.

The values of the model parameters are used to compute the expression of the
reproduction rate, denoted in this paper by $R^{0}$, which represents the
average number of secondary cases produced by one infected individual
introduced into a population of susceptible individuals. This is the crucial
indicator in a transmissible disease. Detailed approaches of the
determination of the reproduction rate in particular models are found in the
literature, see e.g., \cite{Ia-Pu}, \cite{Driessche} and in the references
indicated there.

Thus, the mathematical model we propose here is 
\begin{equation}
S^{\prime }(t)=-\beta _{I}(t)S(t)I(t)-\beta _{A}(t)S(t)A(t)+\xi (t)R(t),%
\mbox{ }  \label{1}
\end{equation}%
\begin{equation}
A^{\prime }(t)=\beta _{I}(t)S(t)I(t)+\beta _{A}(t)S(t)A(t)-(\sigma +\mu
_{A}+l_{A})A(t),\mbox{ }  \label{2}
\end{equation}%
\begin{equation}
I^{\prime }(t)=\sigma A(t)-(\mu _{I}+l_{I})I(t),  \label{3}
\end{equation}%
\begin{equation}
L^{\prime }(t)=l_{A}A(t)+l_{I}I(t)-\mu _{L}L(t),  \label{4}
\end{equation}%
\begin{equation}
R^{\prime }(t)=\mu _{A}A(t)+\mu _{I}I(t)+\mu _{L}L(t)-\xi (t)R(t),  \label{5}
\end{equation}%
for a.a. $t>0,$ with the initial conditions 
\begin{equation}
S(0)=S_{0},\mbox{ }A(0)=A_{0},\mbox{ }I(0)=I_{0},\mbox{ }L(0)=L_{0},\mbox{ }%
R(0)=R_{0}.  \label{6}
\end{equation}%
The model is given here in a normalized form, that is the real number of
individuals at time $t$ in each class is divided to the total population $%
N(t).$

The sum of the individuals in all compartments gives the total population $%
N(t)$ at time $t.$ By (\ref{1})-(\ref{5}) we observe that 
\begin{equation*}
\frac{d}{dt}(S(t)+A(t)+I(t)+L(t)+R(t))=0,
\end{equation*}
which implies that 
\begin{equation}
S(t)+A(t)+I(t)+L(t)+R(t)=S_{0}+A_{0}+I_{0}+L_{0}+R_{0}=N,\mbox{ for all }%
t\geq 0,  \label{8}
\end{equation}%
where the constant $N$ is exactly the total population which, in this model,
remains unchanged at all $t$. It is clear that if each term in these sums is
nonnegative (as representing a fraction of population), then it is bounded
by $N.$ Moreover, since we work with fractions of population, $N$ defined
before equals 1. However, we shall keep it written as $N$ to precisely
indicate where it occurs.

We assume the following conditions for the coefficients of the system: 
\begin{eqnarray}
&&\beta _{I},\mbox{ }\beta _{A},\mbox{ }\xi \in L^{\infty }(0,\infty ),\mbox{
}\beta _{I}(t),\mbox{ }\beta _{A}(t),\mbox{ }\xi (t)\geq 0\mbox{ a.e. }t\geq
0,  \label{7} \\
&&\sigma ,\mbox{ }\mu _{A},\mbox{ }\mu _{I},\mbox{ }\mu _{L},\mbox{ }l_{A},%
\mbox{ }l_{I}\geq 0,  \notag
\end{eqnarray}%
and denote 
\begin{equation}
k_{1}:=\sigma +\mu _{A}+l_{A},\mbox{ }k_{2}:=\mu _{I}+l_{I},\mbox{ }k_{1}>0,%
\mbox{ }k_{2}>0.  \label{11}
\end{equation}

At the end, we make a few comments on some different interpretations of the
model, according to some possible modifications of the coefficients. Thus,
if setting $\xi =0,$ and interpreting $\mu _{A},$ $\mu _{I}$ and $\mu _{L}$
as the mortality of the individuals in the classes $A,$ $I,$ $L,$
respectively, it follows that the class $R$ turns out to correspond to the
extinct population.

If $l_{A}=l_{I}=\mu _{L}=0$ the class of isolated individuals disappear and
so all asymptomatic and symptomatic remain in circulation. Thus, the class $%
L $ is relevant as a control class for the disease containment. A larger
isolation action can be modeled by larger coefficients $l_{A}$ and $l_{I}$.

\subsection{Problem statement}

We assume that two sets of measured values of the isolated and recovered
people, at a time $t=0$ and at a successive time $t=T$, are available.
Namely, it means that we know the nonnegative values $L_{0},$ $R_{0}$ at
time $t=0$ and $L_{T},$ $R_{T},$ at time $t=T$. As specified before, the
number of undetected infected and of the susceptible at these times is not
known, so that a first objective is to estimate $A_{0},$ $I_{0}$, $S_{0}$.
As far as the parameters $\sigma ,\mu _{A},\mu _{I},\mu _{L},l_{A},l_{I},\xi 
$ can be estimated by observations, a direct estimate of the rates $\beta
_{I}$ and $\beta _{A}$ is less obvious. All these justify a study developed
in the present paper, of identifying the rate $\beta _{I}(t)$ and the number
of the undetected infectious individuals $A_{0},$ $I_{0},$ relying on these
available observations for the isolated and recovered people at times $0$
and $T.$ Simultaneously, the number of susceptible $S_{0}$ is identified,
too, because relation (\ref{8}) implies%
\begin{equation}
S_{0}+A_{0}+I_{0}=N_{0}:=N-(L_{0}+R_{0}),  \label{9}
\end{equation}%
whence $S_{0}=N_{0}-(A_{0}+I_{0}).$ Thus, it is sufficient to identify only $%
A_{0}$ and $I_{0}.$

Once the information about the size of the populations $A,$ $I$ and $S$ is
available at time $0$, a prediction about their values at a further time $t$
can be done.

The second objective is to control within a successive time interval the
action of isolating more infected individuals by means of the controllers $%
l_{A}$ and $l_{I}.$ More precisely, the target is to reduce the number of
infected $A$ and $I$ by various actions which can lead to the isolation of
those confirmed, by removing them from circulation and transferring in $L.$
This action is supposed however to be led such that $L$ should not exceed an
upper bound $\widehat{L}.$

These proposed objectives will be expressed by two minimization problems.

\textbf{Problem} $(P_{0}).$ We introduce the cost functional 
\begin{eqnarray}
J(\beta _{I},A_{0},I_{0}) &=&\frac{1}{2}(L(T)-L_{T})^{2}+\frac{1}{2}%
(R(T)-R_{T})^{2}+\frac{\alpha _{1}}{2}\int_{0}^{T}\beta _{I}^{2}(t)dt  \notag
\\
&&+\frac{\alpha _{0}}{2}\left(
A_{0}^{2}+I_{0}^{2}+(N_{0}-A_{0}-I_{0})^{2}\right)   \label{10}
\end{eqnarray}%
and the minimization problem $(P_{0})$ below 
\begin{eqnarray*}
&&\mbox{Minimize }\left\{ J(\beta _{I},A_{0},I_{0});\mbox{ }\beta
_{I}(t)\geq 0\mbox{ a.e. }t\in (0,T),\mbox{ }\right. \mbox{ \ \ \ \ \ \ \ \
\ \ \ \ \ \ \ \ \ \ \ \ \ }(P_{0}) \\
&&\mbox{ \ \ \ }\left. \mbox{ \ \ \ \ \ \ \ \ \ \ \ \ \ }A_{0}\geq 0,\mbox{ }%
I_{0}\geq 0,\mbox{ }A_{0}+I_{0}\leq N_{0},\mbox{ }A(t)\geq 0\mbox{ for all }%
t\in \lbrack 0,T]\right\} 
\end{eqnarray*}%
subject to (\ref{1})-(\ref{6}), (\ref{8}). Here, $\alpha _{0},$ $\alpha _{1}$
are positive constants which may give a larger or smaller weight to the
terms they multiply. The constraint $A_{0}+I_{0}\leq N_{0}$ follows by the
natural assumption that all data are should be nonnegative and so $%
S_{0}=N_{0}-A_{0}-I_{0}\geq 0.$

\textbf{Problem} $(P).$ For the second objective we introduce the cost
functional 
\begin{equation}
J(l_{A},l_{I})=\frac{\alpha _{0}}{2}\int_{0}^{T}A^{2}(t)dt+\frac{\alpha _{0}%
}{2}\int_{0}^{T}I^{2}(t)dt+\frac{\alpha _{1}}{2}(l_{A}^{2}+l_{I}^{2})
\label{10-1}
\end{equation}%
and formulate the optimal control problem%
\begin{equation*}
\mbox{Minimize }J\left\{ (l_{A},l_{I});\mbox{ }l_{A}\in \lbrack 0,1],\mbox{ }%
l_{I}\in \lbrack 0,1],\mbox{ }0\leq L(t)\leq \widehat{L}\right\} \mbox{ \ \
\ \ \ \ \ \ \ \ \ \ \ \ \ \ \ \ \ \ \ }(P)
\end{equation*}%
subject to (\ref{1})-(\ref{6}), (\ref{8}), where $\widehat{L}$ is a fixed
constant, $\widehat{L}>L_{0}.$ The upper bound $\widehat{L}$ for $L$ is
justified by the fact that we try to catch in the class $L$ as much as
possible individuals from the classes $A$ and $I$, but not the total
population. The aim is to detect, by enforcing the testing, at least a part
of the population $A_{0}+I_{0}$ in order to isolate it. It should be said
that the lower bound $L\geq 0$ is not a constraint because this follows from
a property which will be proved for the solution to the state system.

We note that $(P)$ is an optimal control problem with the state constraint $%
L(t)\in \lbrack 0,\widehat{L}]$ which will require a more elaborated
treatment. In fact, for such a problem, the maximum principle (the first
order conditions of optimality) lead to a singular dual backward system.
Problem $(P_{0})$ is much simpler, as we shall see, and this entitles us to
begin our study by approaching first problem $(P)$ and dealing after then
with problem $(P_{0}).$ Thus, we start to solve $(P)$ by assuming that the
values of $A,$ $I,$ $S$ at the time $T$ are calculated after solving $%
(P_{0}) $ and they become the known $A_{0},$ $I_{0},$ $S_{0}$ by resetting
the time at $0.$

We approach problem $(P)$ by an optimal control technique. In Section 2,
after proving the existence and uniqueness of the solution to the state
system (\ref{1})-(\ref{6}), we show that there exists at least a solution $%
(l_{A}^{\ast },l_{I}^{\ast })$ to problem $(P)$ in Proposition 2.2. For this
problem with state restrictions, the optimality conditions cannot be
directly calculated, but via an approximating problem $(P_{\varepsilon })$
indexed along a positive parameter $\varepsilon ,$ which contains
appropriate penalized terms replacing the state constraint. This is
introduced in Section 2.1. The convergence of a sequence of solutions to $%
(P_{\varepsilon })$ precisely to a certain chosen solution to $(P)$ is
proved in Proposition 2.4 and the approximating optimality conditions are
provided in Proposition 2.5. Relying on appropriate estimates for the
solution to the dual system proved in Proposition 2.6, local conditions of
optimality for problem $(P)$ are obtained by passing to the limit in the
approximating ones, in Theorem 2.8. Problem $(P_{0})$ is solved in Section 3
and the optimality conditions are given in Proposition 3.2. An investigation
of the asymptotic stability of the system done in Section 4 finds the
conditions under which the disease can evolve towards an asymptotic
equilibrium state and allows the definition of the reproduction rate. Some
final interpretations in Section 5 complete the paper.

\section{Problem $(P)$}

\setcounter{equation}{0}

We begin with the proof of the well-posedness of the state system. Its
solution will be sometimes denoted by $X:=(S,A,I,L,R).$

By a solution to (\ref{1})-(\ref{6}) on $[0,T]$ we mean an $\mathbb{R}^{5}$%
-valued absolutely continuous function $X=(S,A,I,L,R)$ on $[0,T]$ which
satisfies (\ref{1})-(\ref{6}) a.e. on $(0,T).$

\medskip

\noindent \textbf{Proposition 2.1.} \textit{Let }$T>0$ \textit{and} $%
(S_{0},A_{0},I_{0},L_{0},R_{0})\geq 0.$ \textit{The state system} (\ref{1})-(%
\ref{6})\textit{\ has a unique global solution }$(S,A,I,L,R)\in (W^{1,\infty
}(0,T))^{5}$\textit{. The solution is continuous with respect to the data }$%
(l_{A},l_{I}).$

\medskip

\noindent \textbf{Proof. }In system (\ref{1})-(\ref{5}) we apply the Banach
fixed point theorem, using the set 
\begin{equation*}
\mathcal{M}=\{A\in C([0,T];\mbox{ }0\leq A(t)\leq N\mbox{ for all }t\in
\lbrack 0,T]\}.
\end{equation*}%
Let us pick $a\in \mathcal{M}$ and fix it in the system 
\begin{equation}
S^{\prime }(t)=-\beta _{I}(t)S(t)I(t)-\beta _{A}(t)S(t)a(t)+\xi (t)R(t),%
\mbox{ }  \label{1-1}
\end{equation}%
\begin{equation}
A^{\prime }(t)=\beta _{I}(t)S(t)I(t)+\beta _{A}(t)S(t)a(t)-(\sigma +\mu
_{A}+l_{A})A(t),\mbox{ }  \label{2-1}
\end{equation}%
\begin{equation}
I^{\prime }(t)=\sigma a(t)-(\mu _{I}+l_{I})I(t),  \label{3-1}
\end{equation}%
\begin{equation}
L^{\prime }(t)=l_{A}a(t)+l_{I}I(t)-\mu _{L}L(t),  \label{4-1}
\end{equation}%
\begin{equation}
R^{\prime }(t)=\mu _{A}a(t)+\mu _{I}I(t)+\mu _{L}L(t)-\xi (t)R(t),
\label{5-1}
\end{equation}%
with the initial condition (\ref{6}). We define $\Psi :\mathcal{M}%
\rightarrow C([0,T])$ by $\Psi (a)=A$ where $(S,A,I,L,R)$ is the solution to
(\ref{1-1})-(\ref{5-1}), (\ref{6}) and show that $\Psi (\mathcal{M})\subset 
\mathcal{M}$ and that $\Psi $ is a contraction on $\mathcal{M}.$ By (\ref%
{3-1}) we deduce, using the formula of variation of constants, that 
\begin{equation*}
I(t)=I_{0}e^{-k_{2}t}+\sigma \int_{0}^{t}e^{-k_{2}(t-s)}a(s)ds,\mbox{ }t\in
\lbrack 0,T],
\end{equation*}%
and we have $I\in C([0,T])\cap W^{1,\infty }(0,T).$ Moreover, $I(t)\geq 0$
for all $t\geq 0.$ Applying successively the same formula in (\ref{4-1}), (%
\ref{5-1}), (\ref{1-1}) and (\ref{2-1}) we obtain that $L,R,S,A\in
C([0,T])\cap W^{1,\infty }(0,T)$ and each of them is nonnegative. By (\ref{8}%
) each component is less or equal to $N.$ Thus, $\Psi (\mathcal{M})\subset 
\mathcal{M}.$ It remains to prove that $\Psi $ is a contraction. Let us take
two solutions to (\ref{1-1})-(\ref{5-1}) $(S,A,I,L,R)$ and $(\overline{S},%
\overline{A},\overline{I},\overline{L},\overline{R})$ corresponding to $a$
and $\overline{a}$ respectively, with the same initial condition. By
calculating $(A-\overline{A})(t),$ $(S-\overline{S})(t),...,(R-\overline{R}%
)(t)$ by each corresponding equation we obtain 
\begin{equation*}
\sup_{t\in \lbrack 0,T]}\left\vert (A-\overline{A})(t)\right\vert \leq
C_{1}(T)\sup_{t\in \lbrack 0,T]}\left\vert (a-\overline{a})(t)\right\vert 
\end{equation*}%
where $C_{1}(T)$ is a polynomial in $T$ with coefficients consisting in sums
of the constant systems parameters and the $L^{\infty }$-norms of the time
dependent system parameters. This shows that $\Psi $ is a contraction for
small $T,$ that is we obtain a local solution. Since all solution components
are bounded by $N$ it follows that the solution is global (see e.g., \cite%
{VB-ode}, p. 41, Theorem 2.15).

Let $(l_{A}^{n},l_{I}^{n})_{n}$ be a sequence such that $l_{A}^{n}%
\rightarrow l_{A},$ $l_{I}^{n}\rightarrow l_{I}$ as $n\rightarrow \infty $
and let $X^{n}$ and $X$ be the solutions to (\ref{1})-(\ref{6})
corresponding to these data, respectively. Since the solution $X^{n}\in
(W^{1,\infty }(0,T))^{5}$ it follows that $X^{n}\rightarrow \widetilde{X}$
uniformly in $[0,T]$ and weak* in $(W^{1,\infty }(0,T))^{5}$ and by passing
to the limit in (\ref{1})-(\ref{6}) we deduce that $\widetilde{X}^{\prime
}(t)=\lim_{n\rightarrow \infty }(X^{n})^{\prime }(t),$ whence it follows
that $\widetilde{X}$ is the solution to (\ref{1})-(\ref{6}). \hfill $\square 
$

\medskip

\noindent \textbf{Proposition 2.2. }\textit{Problem }$(P)$\textit{\ has at
least one solution }$(l_{A}^{\ast },l_{I}^{\ast })$.

\medskip

\noindent \textbf{Proof. }It is obvious that an admissible pair exists. For
example, for $l_{A}=l_{I}=0$ we get $L(t)\leq L_{0}<\overline{L}$ and $%
J(l_{A},l_{I})<\infty .$ Let $d:=\inf J(l_{A},l_{I})\geq 0.$ Let us consider
a minimizing sequence $(l_{A}^{n},l_{I}^{n})_{n}$, satisfying the
restrictions in $(P)$. The minimizing sequence also satisfies 
\begin{equation*}
d\leq J(l_{A}^{n},l_{I}^{n})\leq d+\frac{1}{n},\mbox{ for }n\geq 1.
\end{equation*}%
This implies that $l_{A}^{n}\rightarrow l_{A}^{\ast }$ and $%
l_{I}^{n}\rightarrow l_{I}^{\ast }$ as $n\rightarrow \infty $ and $%
l_{A}^{\ast },l_{I}^{\ast }\in \lbrack 0,1].$ We denote by $X^{n}$ the
solution to the state system corresponding to the minimizing sequence$.$ By
Proposition 2.1 this solution exists and belong to $(W^{1,\infty }(0,T))^{5}$
and each component belongs to $[0,N].$ Then, on a subsequence, $%
X^{n}\rightarrow X^{\ast },$ $(X^{n})^{\prime }\rightarrow (X^{\ast
})^{\prime }$ weak* in $(L^{\infty }(0,T))^{5}.$ Therefore, by Arzel\`{a}
theorem it follows that $X^{n}\rightarrow X^{\ast }$ uniformly in $[0,T],$
and so $X^{n}(0)\rightarrow X^{\ast }(0)=X_{0}$ and $X^{n}(T)\rightarrow
X^{\ast }(T).$ According to the last part of Proposition 2.1, it follows
that $X^{\ast }$ is the solution to the state system corresponding to $%
(l_{A}^{\ast },l_{I}^{\ast }).$ Finally, by the weakly lower semicontinuity
of the norms we get $\lim_{n\rightarrow \infty
}J(l_{A}^{n},l_{I}^{n})=J(l_{A}^{\ast },l_{I}^{\ast })$ and so $(l_{A}^{\ast
},l_{I}^{\ast })$ turns out to be optimal in $(P).$\hfill $\square $

\subsection{The approximating problem $(P_{\protect\varepsilon })$}

Let $\varepsilon >0$ and let $(l_{A}^{\ast },l_{I}^{\ast })$ be optimal in $%
(P).$ We introduce the adapted approximating cost functional 
\begin{eqnarray}
J_{\varepsilon }(l_{A},l_{I}) &=&\frac{\alpha _{0}}{2}%
\int_{0}^{T}(A^{2}(t)+I^{2}(t))dt+\frac{\alpha _{1}}{2}(l_{A}^{2}+l_{I}^{2})
\label{13} \\
&&+\frac{\alpha _{2}}{2\varepsilon }\int_{0}^{T}\left( (L(t)-\widehat{L}%
)^{+}\right) ^{2}dt+\frac{1}{2}(l_{A}-l_{A}^{\ast })^{2}+\frac{1}{2}%
(l_{I}-l_{I}^{\ast })^{2}  \notag
\end{eqnarray}%
and study the following approximating problem $(P_{\varepsilon }),$ 
\begin{equation*}
\begin{tabular}{ll}
Minimize $\left\{ J_{\varepsilon }(l_{A},l_{I});\mbox{ }l_{A}\in \lbrack
0,1],\mbox{ }l_{I}\in \lbrack 0,1]\right\} $ & $\mbox{ \ \ \ \ \ \ \ \ \ \ \
\ \ \ \ \ \ \ \ \ \ \ \ \ \ \ \ \ \ \ \ \ \ \ }(P_{\varepsilon })$%
\end{tabular}%
\end{equation*}%
subject to (\ref{1})-(\ref{6}), (\ref{8}).

We observe that the state constraint in $(P)$ is replaced here by the
penalization of the $L^{2}$-norm of the positive part of $(L(t)-\widehat{L}%
), $ where $\alpha _{2}>0.$ The last two penalization terms in (\ref{13})
ensure the convergence of the approximating solution to the chosen optimal
controller $(l_{A}^{\ast },l_{I}^{\ast })$ in $(P)$.

\medskip

\noindent \textbf{Proposition 2.3. }\textit{Problem }$(P_{\varepsilon })$%
\textit{\ has at least one solution, }$(l_{A,\varepsilon }^{\ast
},l_{I,\varepsilon }^{\ast })$ \textit{with the corresponding state} $%
X_{\varepsilon }^{\ast }\in (W^{1,\infty }(0,T))^{5},$ \textit{having the
components in} $[0,N].$

\medskip \medskip

\noindent \textbf{Proof. }First of all we see that there is at least an
admissible triplet, let it be $(l_{A}^{\ast \ast },l_{I}^{\ast \ast })$ an
optimal one in $(P)$ with the corresponding global state $X^{\ast \ast }$,
belonging to $(W^{1,\infty }(0,T))^{5}$ and $L^{\ast \ast }\leq \widehat{L}.$
Then, $J_{\varepsilon }(l_{A}^{\ast \ast },l_{I}^{\ast \ast })<\infty .$
Hence, the admissible set in $(P_{\varepsilon })$ is not empty and since $%
J_{\varepsilon }(l_{A},l_{I})\geq 0$, there exists $d_{\varepsilon }:=\inf
J_{\varepsilon }(l_{A},l_{I})\geq 0.$ We take a minimizing sequence $%
(l_{A,\varepsilon }^{n},l_{I,\varepsilon }^{n})$ with the corresponding
solution $X_{\varepsilon }^{n}$ to the state system in the class of global
solutions, with the components in the interval $[0,N]$. Recall that the
initial condition $X_{0}^{n}=X_{0}$ is nonnegative. We have $d_{\varepsilon
}\leq J_{\varepsilon }(l_{A,\varepsilon }^{n},l_{I,\varepsilon }^{n})\leq d+%
\frac{1}{n},$ for $n\geq 1.$ Then, $l_{A,\varepsilon }^{n}\rightarrow
l_{A,\varepsilon }^{\ast },$ $l_{I,\varepsilon }^{n}\rightarrow
l_{I,\varepsilon }^{\ast }$ and $l_{A,\varepsilon }^{\ast },l_{I,\varepsilon
}^{\ast }\in \lbrack 0,1].$ By (\ref{13}) it follows that $(L_{\varepsilon
}^{n}-\widehat{L})^{+},$ $(A_{\varepsilon }^{n})_{n},$ $(I_{\varepsilon
}^{n})_{n}$ are bounded in $L^{2}(0,T)$ and so by (\ref{4}), (\ref{5}), (\ref%
{1})-(\ref{3}) we get that $(L_{\varepsilon }^{n})_{n}$, $(R_{\varepsilon
}^{n})_{n},$ $(S_{\varepsilon }^{n})_{n},$ $(A_{\varepsilon }^{n})_{n},$ $%
(I_{\varepsilon }^{n})_{n}$ are bounded in $W^{1,\infty }(0,T)$. We infer
that $X_{\varepsilon }^{n}\rightarrow X_{\varepsilon }^{\ast }$ uniformly in 
$[0,T]$. The limit is bounded and satisfies the state system corresponding
to $(l_{A}^{\ast },l_{I}^{\ast }).$ In addition we note that since $%
L\rightarrow L^{+}$ is continuous, we have $(L_{\varepsilon }^{n}-\widehat{L}%
)^{+}\rightarrow (L_{\varepsilon }^{\ast }-\widehat{L})^{+}$ uniformly in $%
[0,T]$ and so, $\lim_{n\rightarrow \infty }J_{\varepsilon }(l_{A,\varepsilon
}^{n},l_{I,\varepsilon }^{n})=J_{\varepsilon }(l_{A,\varepsilon }^{\ast
},l_{I,\varepsilon }^{\ast }).$ Moreover, all components are nonnegative and
less or equal to $N.$ All these prove that $(l_{A,\varepsilon }^{\ast
},l_{I,\varepsilon }^{\ast })$ is optimal in $(P_{\varepsilon }).$\hfill $%
\square $

\medskip

\noindent \textbf{Proposition 2.4. }\textit{Let }$\{(l_{A}^{\ast
},l_{I}^{\ast }),$ $X^{\ast }\}$\textit{\ and }$\mathit{\{}(l_{A,\varepsilon
}^{\ast },l_{I,\varepsilon }^{\ast }),X_{\varepsilon }^{\ast
}\}_{\varepsilon }$\textit{\ be optimal in }$(P)$ \textit{and in }$%
(P_{\varepsilon }),$ \textit{respectively. Then,} 
\begin{equation}
l_{A,\varepsilon }^{\ast }\rightarrow l_{A}^{\ast },\mbox{ }l_{I,\varepsilon
}^{\ast }\rightarrow l_{I}^{\ast }\mbox{, \textit{as} }\varepsilon
\rightarrow 0,  \label{13-0}
\end{equation}%
\begin{equation}
X_{\varepsilon }\rightarrow X^{\ast }\mbox{ \textit{weak* in} }(W^{1,\infty
}(0,T))^{5}\mbox{ \textit{and uniformly in} }[0,T],\mbox{ , \textit{as} }%
\varepsilon \rightarrow 0.  \label{13-1}
\end{equation}

\medskip

\noindent \textbf{Pro}of. If $(l_{A,\varepsilon }^{\ast },l_{I,\varepsilon
}^{\ast })$ is optimal, then $J_{\varepsilon }(l_{A,\varepsilon }^{\ast
},l_{I,\varepsilon }^{\ast })\leq J_{\varepsilon }(l_{A},l_{I}),$ for all $%
(l_{A},l_{I})$ satisfying the constraints in $(P_{\varepsilon }).$ In
particular, we can set $(l_{A},l_{I})=(l_{A}^{\ast },l_{I}^{\ast })$ which
is the optimal triplet chosen in $(P)$ and then the previous inequality
becomes 
\begin{eqnarray}
&&J_{\varepsilon }(l_{A,\varepsilon }^{\ast },l_{I,\varepsilon }^{\ast })=%
\frac{\alpha _{0}}{2}\int_{0}^{T}((A_{\varepsilon }^{\ast
}(t))^{2}+(I_{\varepsilon }^{\ast }(t))^{2})dt+\frac{\alpha _{1}}{2}%
(l_{A,\varepsilon }^{2}+l_{I,\varepsilon }^{2})  \label{14} \\
&&+\frac{\alpha _{2}}{2\varepsilon }\int_{0}^{T}\left( (L_{\varepsilon
}^{\ast }(t)-\widehat{L})^{+}\right) ^{2}(t)dt+\frac{1}{2}(l_{A,\varepsilon
}^{\ast }-l_{A}^{\ast })^{2}+\frac{1}{2}(l_{I,\varepsilon }^{\ast
}-l_{I}^{\ast })^{2}  \notag \\
&\leq &J_{\varepsilon }(l_{A}^{\ast },l_{I}^{\ast })=\frac{\alpha _{0}}{2}%
\int_{0}^{T}((A^{\ast }(t))^{2}+(I^{\ast }(t))^{2})dt+\frac{\alpha _{1}}{2}%
(l_{A}^{2}+l_{I}^{2}),  \notag
\end{eqnarray}%
because $L^{\ast }\leq \widehat{L}$ in $(P).$ Thus, the left-hand side is
bounded and we have $l_{A,\varepsilon }^{\ast }\rightarrow l_{A}^{\ast },$ $%
l_{I,\varepsilon }^{\ast }\rightarrow l_{I}^{\ast }$, as $\varepsilon
\rightarrow 0.$ By Proposition 2.1, we infer that $X_{\varepsilon
}\rightarrow X^{\ast }$ weak* in $(W^{1,\infty }(0,T))^{5}$ and uniformly in 
$[0,T].$ By (\ref{14}) we can write 
\begin{equation*}
\frac{\alpha _{2}}{2\varepsilon }\int_{0}^{T}\left( (L_{\varepsilon }^{\ast
}(t)-\widehat{L})^{+}\right) ^{2}(t)dt\leq \mbox{ constant,}
\end{equation*}%
which implies that $\int_{0}^{T}\left( (L_{\varepsilon }^{\ast }(t)-\widehat{%
L})^{+}\right) ^{2}(t)dt\rightarrow 0$ and so, $(L_{\varepsilon }^{\ast }(t)-%
\widehat{L})^{+}\rightarrow 0$ strongly in $L^{2}(0,T).$ On the other hand, $%
(L_{\varepsilon }^{\ast }-\widehat{L})^{+}(t)\rightarrow (L^{\ast }-\widehat{%
L})^{+}(t)$ for all $t\in \lbrack 0,T],$ so that $L^{\ast }(t)\leq \widehat{L%
},$ which ends the proof$.$\hfill $\square $

\subsection{The approximating optimality conditions}

Let $\lambda >0$ and set the variations%
\begin{eqnarray*}
l_{A,\varepsilon }^{\lambda } &=&l_{A,\varepsilon }^{\ast }+\lambda \omega
_{A},\mbox{ }\omega _{A}=\widetilde{l_{A}}-l_{A,\varepsilon }^{\ast },\mbox{ 
}\widetilde{l_{A}}\in \lbrack 0,1], \\
l_{I,\varepsilon }^{\lambda } &=&l_{I,\varepsilon }^{\ast }+\lambda \omega
_{I},\mbox{ }\omega _{I}=\widetilde{l_{I}}-l_{I,\varepsilon }^{\ast },\mbox{ 
}\widetilde{l_{I}}\in \lbrack 0,1].\mbox{ }
\end{eqnarray*}%
Let us denote $x_{\varepsilon }^{\lambda }:=\frac{X_{\varepsilon
}^{l_{A,\varepsilon }^{\lambda },l_{I,\varepsilon }^{\lambda
}}-X_{\varepsilon }^{\ast }}{\lambda },$ where $X_{\varepsilon
}^{l_{A,\varepsilon }^{\lambda },l_{I,\varepsilon }^{\lambda }}$ is the
solution to (\ref{1})-(\ref{6}) corresponding to $(l_{A,\varepsilon
}^{\lambda },l_{I,\varepsilon }^{\lambda }),$ satisfying (\ref{8}) and $%
X_{\varepsilon }^{\ast }$ is the optimal state corresponding to $%
(l_{A,\varepsilon }^{\ast },l_{I,\varepsilon }^{\ast }).$

We introduce the linearized system for problem $(P)$%
\begin{eqnarray}
s^{\prime } &=&-k_{0,\varepsilon }s-\beta _{I}S_{\varepsilon }^{\ast
}i-\beta _{A}S_{\varepsilon }^{\ast }a+\xi r,  \notag \\
a^{\prime } &=&k_{0,\varepsilon }s+\beta _{I}^{\ast }S_{\varepsilon }^{\ast
}i+k_{3,\varepsilon }a-\omega _{A}A_{\varepsilon }^{\ast },  \notag \\
i^{\prime } &=&\sigma a-k_{2,\varepsilon }i-\omega _{I}I_{\varepsilon
}^{\ast },  \label{15} \\
l^{\prime } &=&l_{A,\varepsilon }^{\ast }a+l_{I,\varepsilon }^{\ast }i-\mu
_{L}l+\omega _{A}A_{\varepsilon }^{\ast }+\omega _{I}I_{\varepsilon }^{\ast
},  \notag \\
r^{\prime } &=&\mu _{A}a+\mu _{I}i+\mu _{L}l-\xi r,  \notag
\end{eqnarray}%
for a.a. $t>0,$ with the initial condition 
\begin{equation}
s(0)=0,\mbox{ }a(0)=0,\mbox{ }i(0)=0,\mbox{ }l(0)=0,\mbox{ }r(0)=0,
\label{15-0}
\end{equation}%
where 
\begin{eqnarray}
k_{0,\varepsilon } &=&\beta _{A}A_{\varepsilon }^{\ast }+\beta
_{I}I_{\varepsilon }^{\ast },\mbox{ }k_{1,\varepsilon }=\sigma +\mu
_{A}+l_{A,\varepsilon }^{\ast },\mbox{ }  \label{16} \\
k_{2,\varepsilon } &=&\mu _{I}+l_{I,\varepsilon }^{\ast },\mbox{ }%
k_{3,\varepsilon }=\beta _{A}S_{\varepsilon }^{\ast }-k_{1,\varepsilon }. 
\notag
\end{eqnarray}%
First, by the known results for linear systems we infer that (\ref{15}) has
a unique global solution $x=(s,a,i,l,r)\in (W^{1,\infty }(0,T))^{5}.$

By a direct calculation, using the continuity with respect to the data of
the solution to the state system it can be easily proved that $%
x_{\varepsilon }^{\lambda }\rightarrow x_{\varepsilon }:=(s,a,i,l,r)$
strongly in $C([0,T])$, as $\lambda \rightarrow 0,$ so that (\ref{15})
stands for the system in variations.

\noindent We introduce the backward dual system for the variables $%
(p_{\varepsilon },q_{\varepsilon },d_{\varepsilon },e_{\varepsilon
},f_{\varepsilon })$ as 
\begin{equation}
p_{\varepsilon }^{\prime }-k_{0,\varepsilon }p_{\varepsilon
}+k_{0,\varepsilon }q_{\varepsilon }=0,  \label{17}
\end{equation}%
\begin{equation}
q_{\varepsilon }^{\prime }-\beta _{A}S_{\varepsilon }^{\ast }p_{\varepsilon
}+k_{3,\varepsilon }q_{\varepsilon }+\sigma d_{\varepsilon
}+l_{A,\varepsilon }^{\ast }e_{\varepsilon }+\mu _{A}f_{\varepsilon
}=-\alpha _{0}A_{\varepsilon }^{\ast },  \label{18}
\end{equation}%
\begin{equation}
d_{\varepsilon }^{\prime }-\beta _{I}S_{\varepsilon }^{\ast }p_{\varepsilon
}+\beta _{I}S_{\varepsilon }^{\ast }q_{\varepsilon }-k_{2,\varepsilon
}d_{\varepsilon }+l_{I,\varepsilon }^{\ast }e_{\varepsilon }+\mu
_{I}f_{\varepsilon }=-\alpha _{0}I_{\varepsilon }^{\ast },  \label{19}
\end{equation}%
\begin{equation}
e_{\varepsilon }^{\prime }-\mu _{L}e_{\varepsilon }+\mu _{L}f_{\varepsilon }+%
\frac{\alpha _{2}}{\varepsilon }(L_{\varepsilon }^{\ast }-\widehat{L})^{+}=0,
\label{20}
\end{equation}%
\begin{equation}
f_{\varepsilon }^{\prime }+\xi p_{\varepsilon }-\xi f_{\varepsilon }=0,
\label{21}
\end{equation}%
for a.a. $t>0,$ with the final conditions 
\begin{equation}
p_{\varepsilon }(T)=0,\mbox{ }q_{\varepsilon }(T)=0,\mbox{ }d_{\varepsilon
}(T)=0,\mbox{ }e_{\varepsilon }(T)=0,\mbox{ }f_{\varepsilon }(T)=0.
\label{22}
\end{equation}%
The linear system (\ref{17})-(\ref{22}) has, for each $\varepsilon >0,$ a
unique global solution $(p_{\varepsilon },q_{\varepsilon },d_{\varepsilon
},e_{\varepsilon },f_{\varepsilon })\in (W^{1,\infty }(0,T))^{5}.$

Let $N_{[0,1]}(z)$ be the normal cone to the set $[0,1]$, 
\begin{equation*}
N_{[0,1]}(z)=\left\{ 
\begin{array}{l}
\mathbb{R}^{-},\mbox{ if }z=0 \\ 
0,\mbox{ \ \ \ if }z\in (0,1), \\ 
\mathbb{R}^{+},\mbox{ if }z=1.%
\end{array}%
\right.
\end{equation*}

We recall that the projection of a point $z$ on a set $K\subset \mathbb{R}%
^{d},$ $d\in \mathbb{N},$ is defined by $P_{K}(z):=(I_{d}+\kappa \partial
I_{K})^{-1}(z)$, for all $k>0,$ $I_{d}$ being the identity operator.

\medskip

\noindent \textbf{Proposition 2.5. }\textit{Let }$(l_{A,\varepsilon }^{\ast
},l_{I,\varepsilon }^{\ast })$\textit{\ be optimal in }$(P_{\varepsilon })$ 
\textit{with the state }$X_{\varepsilon }^{\ast }.$\textit{\ Then, }%
\begin{eqnarray}
l_{A,\varepsilon }^{\ast } &=&P_{[0,1]}\left( \frac{1}{\alpha _{1}+1}\left(
\int_{0}^{T}A_{\varepsilon }^{\ast }(q_{\varepsilon }-e_{\varepsilon
})dt+l_{A}^{\ast }\right) \right) ,  \label{23} \\
l_{I,\varepsilon }^{\ast } &=&P_{[0,1]}\left( \frac{1}{\alpha _{1}+1}\left(
\int_{0}^{T}I_{\varepsilon }^{\ast }(d_{\varepsilon }-e_{\varepsilon
})dt+l_{I}^{\ast }\right) \right) ,  \notag
\end{eqnarray}%
\textit{where} $(p_{\varepsilon },q_{\varepsilon },d_{\varepsilon
},e_{\varepsilon },f_{\varepsilon })$\textit{\ is the solution to the
backward dual system} (\ref{17})-(\ref{22}).

\medskip

\noindent \textbf{Proof. }Let us multiply the equations for $s,a,i,l,r$ in (%
\ref{15}) by $p_{\varepsilon },q_{\varepsilon },d_{\varepsilon
},e_{\varepsilon },f_{\varepsilon },$ respectively and integrate over $(0,T).
$ By integrating by parts and taking into account the equations in the dual
system and the initial conditions in the system in variations, we obtain%
\begin{eqnarray}
&&\alpha _{0}\int_{0}^{T}(A_{\varepsilon }^{\ast }a+I_{\varepsilon }^{\ast
}i)dt+\frac{\alpha _{2}}{\varepsilon }\int_{0}^{T}(L_{\varepsilon }^{\ast }-%
\widehat{L})^{+}ldt  \label{25} \\
&=&\int_{0}^{T}(\omega _{A}A_{\varepsilon }^{\ast }(e_{\varepsilon
}-q_{\varepsilon })+\omega _{I}I_{\varepsilon }^{\ast }(e_{\varepsilon
}-d_{\varepsilon }))dt.  \notag
\end{eqnarray}%
On the other hand, for $(l_{A,\varepsilon }^{\ast },l_{I,\varepsilon }^{\ast
})$ optimal in $(P_{\varepsilon })$ we can write 
\begin{equation*}
J_{\varepsilon }(l_{A,\varepsilon }^{\lambda },l_{I,\varepsilon }^{\lambda
})\geq J_{\varepsilon }(l_{A,\varepsilon }^{\ast },l_{I,\varepsilon }^{\ast
}).
\end{equation*}%
By replacing the expression of the cost functional $J_{\varepsilon },$
performing some algebra, dividing by $\lambda $ and passing to the limit as $%
\lambda \rightarrow 0$ we obtain 
\begin{eqnarray}
&&\alpha _{0}\int_{0}^{T}(A_{\varepsilon }^{\ast }a+I_{\varepsilon }^{\ast
}i)dt+\alpha _{1}(l_{A,\varepsilon }^{\ast }\omega _{A}+l_{I,\varepsilon
}^{\ast }\omega _{I})  \label{26} \\
&&+\frac{\alpha _{2}}{\varepsilon }\int_{0}^{T}(L_{\varepsilon }^{\ast }-%
\widehat{L})^{+}ldt+(l_{A,\varepsilon }^{\ast }-l_{A}^{\ast })\omega
_{A}+(l_{I,\varepsilon }^{\ast }-l_{I}^{\ast })\omega _{I}\geq 0.  \notag
\end{eqnarray}%
By comparison with (\ref{25}) we deduce 
\begin{eqnarray}
&&\int_{0}^{T}\omega _{A}A_{\varepsilon }^{\ast }(e_{\varepsilon
}-q_{\varepsilon })dt+\left( (\alpha _{1}+1)l_{A,\varepsilon }^{\ast
}-l_{A}^{\ast }\right) \omega _{A}  \label{27} \\
&&+\int_{0}^{T}\omega _{I}I_{\varepsilon }^{\ast }(e_{\varepsilon
}-d_{\varepsilon })dt+\left( (\alpha _{1}+1)l_{I,\varepsilon }^{\ast
}-l_{I}^{\ast }\right) \omega _{I}\geq 0.  \notag
\end{eqnarray}%
Recall the setting of $\omega _{A}$ and $\omega _{I}$ and choose, in
particular, $\widetilde{l_{A}}=l_{A,\varepsilon }^{\ast },$ meaning that we
keep $l_{A,\varepsilon }^{\ast }$ fixed and give a variation only to $%
l_{I,\varepsilon }^{\ast }.$ Then, (\ref{27}) yields 
\begin{equation*}
\left( -\int_{0}^{T}I_{\varepsilon }^{\ast }(e_{\varepsilon }-d_{\varepsilon
})dt-(\alpha _{1}+1)l_{I,\varepsilon }^{\ast }+l_{I}^{\ast }\right)
(l_{I,\varepsilon }^{\ast }-\widetilde{l_{I}})\geq 0
\end{equation*}%
for all $\widetilde{l_{I}}\in \lbrack 0,1],$ which implies%
\begin{equation}
-\int_{0}^{T}I_{\varepsilon }^{\ast }(e_{\varepsilon }-d_{\varepsilon
})dt-(\alpha _{1}+1)l_{I,\varepsilon }^{\ast }+l_{I}^{\ast }\in
N_{[0,1]}(l_{I,\varepsilon }^{\ast }),  \label{28}
\end{equation}%
whence we have the second relation in (\ref{23}).

Then, let us set $\widetilde{l_{I}}=l_{I,\varepsilon }^{\ast }.$ By (\ref{27}%
) we obtain that 
\begin{equation}
-\int_{0}^{T}A_{\varepsilon }^{\ast }(e_{\varepsilon }-q_{\varepsilon
})dt-(\alpha _{1}+1)l_{A,\varepsilon }^{\ast }+l_{A}^{\ast }\in
N_{[0,1]}(l_{A,\varepsilon }^{\ast })  \label{29-0}
\end{equation}%
which yields the first relation in (\ref{23}).\hfill $\square $

\subsection{The optimality conditions for problem $(P)$}

We begin by proving the boundedness of the solution to the dual system.

\medskip

\noindent \textbf{Proposition 2.6. }\textit{There exists\ }$T_{loc},$\textit{%
\ such that for }$T<T_{loc}$ \textit{we have} 
\begin{equation}
\int_{0}^{T}\left\vert e_{\varepsilon }^{\prime }(t)\right\vert dt\leq C,%
\mbox{ }  \label{29-2}
\end{equation}%
\begin{equation}
\left\vert p_{\varepsilon }(t)\right\vert +\left\vert q_{\varepsilon
}(t)\right\vert +\left\vert d_{\varepsilon }(t)\right\vert +\left\vert
f_{\varepsilon }(t)\right\vert \leq C,\mbox{ \textit{for} }t\in (0,T),
\label{29-3}
\end{equation}%
\textit{independently of} $\varepsilon .$

\medskip

\noindent \textbf{Proof. }Let us choose\textbf{\ }$y_{1}<\widehat{L},$ $\rho
>0$ and $\theta \in \mathbb{R},$ with $\left\vert \theta \right\vert =1$
such that $y_{1}+\rho \theta \leq \widehat{L},$ where $L_{0}<\widehat{L}.$
Then, we have $y_{1}\pm \rho \leq \widehat{L}$ and this holds if $\rho \leq 
\widehat{L}-y_{1}.$ The value $y_{1}$ is at our free choice and we choose $%
0<y_{1}<L_{0}$. We multiply (\ref{20}) by $(L_{\varepsilon }^{\ast
}-(y_{1}+\rho \theta ))$ and integrate over $(0,T).$ We have%
\begin{eqnarray}
&&\int_{0}^{T}(L_{\varepsilon }^{\ast }-y_{1})e_{\varepsilon }^{\prime
}dt-\int_{0}^{T}\rho \theta e_{\varepsilon }^{\prime }dt+\frac{\alpha _{2}}{%
\varepsilon }\int_{0}^{T}(L_{\varepsilon }^{\ast }-\widehat{L}%
)^{+}(L_{\varepsilon }^{\ast }-y_{1}-\rho \theta )dt  \notag \\
&=&\int_{0}^{T}(\mu _{L}e_{\varepsilon }-\mu _{L}f_{\varepsilon
})(L_{\varepsilon }^{\ast }-y_{1}-\rho \theta )dt.  \label{30}
\end{eqnarray}%
Let us define $\varphi (y)=\frac{1}{2}(y^{+})^{2},$ $\varphi :[0,\infty
)\rightarrow \lbrack 0,\infty ).$ Its subdifferential $\partial \varphi
(y)=y^{+}$ and $\varphi (y)-\varphi (z)\leq \partial \varphi (y)(y-z)$ for
all $z\leq \widehat{L}.$ In the second term on the left-hand side of (\ref%
{30}) we set 
\begin{equation*}
\theta :=-\frac{e_{\varepsilon }^{\prime }}{\left\vert e_{\varepsilon
}^{\prime }\right\vert }1_{\{t;e_{\varepsilon }^{\prime }(t)\neq 0\}},
\end{equation*}%
where $1_{\{t;e_{\varepsilon }^{\prime }(t)\neq 0\}}$ is the characteristic
function of the set indicated as subscript. Using the relation for the
subdifferential of $\varphi $ and integrating by parts the first term on the
left-hand side of (\ref{30}) we have%
\begin{eqnarray*}
&&(L_{\varepsilon }^{\ast }(T)-y_{1})e_{\varepsilon }(T)-(L_{\varepsilon
}^{\ast }(0)-y_{1})e_{\varepsilon }(0)-\int_{0}^{T}(L_{\varepsilon }^{\ast
})^{\prime }e_{\varepsilon }dt+\rho \int_{0}^{T}\left\vert e_{\varepsilon
}^{\prime }(t)\right\vert dt \\
&&+\frac{\alpha _{2}}{\varepsilon }\int_{0}^{T}((L_{\varepsilon }^{\ast }-%
\widehat{L})^{+})^{2}dt-\frac{\alpha _{2}}{\varepsilon }%
\int_{0}^{T}(((y_{1}+\rho \theta )-\widehat{L})^{+})^{2}dt \\
&\leq &\int_{0}^{T}\mu _{L}e_{\varepsilon }-\mu _{L}f_{\varepsilon
})L_{\varepsilon }^{\ast }dt-\int_{0}^{T}(\mu _{L}e_{\varepsilon }-\mu
_{L}f_{\varepsilon })(y_{1}+\rho \theta )dt.
\end{eqnarray*}%
Since $y_{1}+\rho \theta \leq \widehat{L}$, $e_{\varepsilon }(T)=0,$ using
eq. (\ref{4}) and making some rearrangements we are led to the relation%
\begin{eqnarray*}
&&\rho \int_{0}^{T}\left\vert e_{\varepsilon }^{\prime }(t)\right\vert
dt\leq \left\vert L_{0}-y_{1}\right\vert \left\vert e_{\varepsilon
}(0)\right\vert  \\
&&+\int_{0}^{T}(e_{\varepsilon }(l_{A,\varepsilon }^{\ast }A_{\varepsilon
}^{\ast }+l_{I,\varepsilon }^{\ast }I_{\varepsilon }^{\ast })-\mu
_{L}f_{\varepsilon }L_{\varepsilon }^{\ast })dt-\int_{0}^{T}(\mu
_{L}e_{\varepsilon }-\mu _{L}f_{\varepsilon })(y_{1}+\rho \theta )dt \\
&\leq &\int_{0}^{T}(\left\vert e_{\varepsilon }\right\vert +\left\vert
f_{\varepsilon })\right\vert )\left( l_{A,\varepsilon }^{\ast }\left\vert
A_{\varepsilon }^{\ast }\right\vert +l_{I,\varepsilon }^{\ast }\left\vert
I_{\varepsilon }^{\ast }\right\vert +\mu _{L}\left\vert L_{\varepsilon
}^{\ast }\right\vert +\mu _{L}y_{1}+\mu _{L}\rho \right) dt.
\end{eqnarray*}%
Denoting 
\begin{equation}
E_{\varepsilon }:=\int_{0}^{T}\left\vert e_{\varepsilon }^{\prime
}(t)\right\vert dt,  \label{30-0}
\end{equation}%
\begin{equation*}
F_{0}=\left\vert L_{0}-y_{1}\right\vert ,\mbox{ }F_{1,\varepsilon
}=l_{A,\varepsilon }^{\ast }\left\Vert A_{\varepsilon }^{\ast }\right\Vert
_{\infty }+l_{I,\varepsilon }^{\ast }\left\Vert I_{\varepsilon }^{\ast
}\right\Vert _{\infty }+\mu _{L}\left\Vert L_{\varepsilon }^{\ast
}\right\Vert _{\infty }+\mu _{L}y_{1}
\end{equation*}%
where $\left\Vert \cdot \right\Vert _{\infty }=\left\Vert \cdot \right\Vert
_{L^{\infty }(0,T)},$ we have 
\begin{equation}
\rho \int_{0}^{T}\left\vert e_{\varepsilon }^{\prime }(t)\right\vert dt\leq
F_{0}\left\vert e_{\varepsilon }(0)\right\vert +\int_{0}^{T}(\left\vert
e_{\varepsilon }\right\vert +\left\vert f_{\varepsilon })\right\vert
)(F_{1,\varepsilon }+\mu _{L}\rho )dt.  \label{31}
\end{equation}%
Now, since $e_{\varepsilon }(T)=0,$ we note that 
\begin{equation}
\left\vert e_{\varepsilon }(0)\right\vert =\left\vert
\int_{0}^{T}e_{\varepsilon }^{\prime }(t)dt\right\vert \leq E_{\varepsilon },%
\mbox{ }\left\vert e_{\varepsilon }(t)\right\vert \leq 2E_{\varepsilon }%
\mbox{ for all }t\in \lbrack 0,T].  \label{38}
\end{equation}%
We multiply eqs. (\ref{17})-(\ref{19}), (\ref{21}) by $p_{\varepsilon
},q_{\varepsilon },d_{\varepsilon },f_{\varepsilon }$ respectively,
integrate over $(0,t)$ and sum up. We get 
\begin{eqnarray*}
\frac{1}{2}\mathcal{S}_{\varepsilon }^{2}(t) &\leq &\frac{1}{2}\alpha
_{0}^{2}\left( \left\Vert A_{\varepsilon }^{\ast }\right\Vert _{\infty
}^{2}+\left\Vert I_{\varepsilon }^{\ast }\right\Vert _{\infty }^{2}\right)
T+\int_{0}^{t}\digamma _{\varepsilon }\mathcal{S}_{\varepsilon }^{2}(\tau
)d\tau  \\
&&+l_{A,\varepsilon }^{\ast }\int_{0}^{t}\left\vert e_{\varepsilon
}\right\vert \left\vert q_{\varepsilon }\right\vert d\tau +l_{I,\varepsilon
}^{\ast }\int_{0}^{t}\left\vert e_{\varepsilon }\right\vert \left\vert
d_{\varepsilon }\right\vert d\tau ,
\end{eqnarray*}%
whence 
\begin{equation*}
\mathcal{S}_{\varepsilon }^{2}(t)\leq \alpha _{0}^{2}\left( \left\Vert
A_{\varepsilon }^{\ast }\right\Vert _{\infty }^{2}+\left\Vert I_{\varepsilon
}^{\ast }\right\Vert _{\infty }^{2}\right) T+\int_{0}^{t}2G_{\varepsilon }%
\mathcal{S}_{\varepsilon }^{2}(\tau )d\tau +8E_{\varepsilon }^{2}T,
\end{equation*}%
where $\mathcal{S}_{\varepsilon }^{2}(t)=\left\vert p_{\varepsilon
}(t)\right\vert ^{2}+\left\vert q_{\varepsilon }(t)\right\vert
^{2}+\left\vert d_{\varepsilon }(t)\right\vert ^{2}+\left\vert
f_{\varepsilon }(t)\right\vert ^{2}.$ Here, $\digamma _{\varepsilon }$ and $%
G_{\varepsilon }$ consist in sums of the constant coefficients of the
equations in the dual system plus the $L^{\infty }$-norms of the time
dependent coefficients. These sums also include $l_{A,\varepsilon }^{\ast 2}$
and $l_{I,\varepsilon }^{\ast 2}$. By the Gronwall's lemma, the previous
inequality yields the estimate 
\begin{equation}
\left\vert p_{\varepsilon }(t)\right\vert +\left\vert q_{\varepsilon
}(t)\right\vert +\left\vert d_{\varepsilon }(t)\right\vert +\left\vert
f_{\varepsilon }(t)\right\vert \leq \alpha _{0}(\left\Vert A_{\varepsilon
}^{\ast }\right\Vert _{\infty }+\left\Vert I_{\varepsilon }^{\ast
}\right\Vert _{\infty }+4E_{\varepsilon })\sqrt{T}e^{G_{\varepsilon }T},
\label{32}
\end{equation}%
the right-hand side being bounded independently of $\varepsilon ,$ since $%
A_{\varepsilon }^{\ast },$ $I_{\varepsilon }^{\ast }$ tend uniformly to $%
A^{\ast },$ $I^{\ast },$ $l_{A,\varepsilon }^{\ast }\rightarrow l_{A}^{\ast
},$ $l_{I,\varepsilon }^{\ast }\rightarrow l_{I}^{\ast },$ $G_{\varepsilon
}\rightarrow G$ which is constant$.$

We go back to (\ref{31}) and using (\ref{32}) we write 
\begin{eqnarray*}
\rho E_{\varepsilon } &\leq &F_{0}E_{\varepsilon }+(\left\vert
e_{\varepsilon }\right\vert +\left\vert f_{\varepsilon })\right\vert
)F_{1,\varepsilon }T+\mu _{L}\rho (\left\vert e_{\varepsilon }\right\vert
+\left\vert f_{\varepsilon })\right\vert )T \\
&\leq &(F_{0}+4\alpha _{0}F_{1,\varepsilon }T\sqrt{T}e^{G_{\varepsilon
}T})E_{\varepsilon }+4\mu _{L}\rho \alpha _{0}T\sqrt{T}e^{G_{\varepsilon
}T}E_{\varepsilon } \\
&&+\alpha _{0}F_{2,\varepsilon }(F_{1,\varepsilon }+\mu _{L}\rho )T\sqrt{T}%
e^{G_{\varepsilon }T},
\end{eqnarray*}%
where $F_{2.\varepsilon }:=\left\Vert A_{\varepsilon }^{\ast }\right\Vert
_{\infty }+\left\Vert I_{\varepsilon }^{\ast }\right\Vert _{\infty }.$ This
implies 
\begin{eqnarray}
&&E_{\varepsilon }\left( \rho (1-4\mu _{L}\alpha _{0}e^{G_{\varepsilon }T}T%
\sqrt{T})-(F_{0}+4\alpha _{0}F_{1,\varepsilon }T\sqrt{T}e^{G_{\varepsilon
}T})\right)  \label{41} \\
&\leq &\alpha _{0}F_{2,\varepsilon }(F_{1,\varepsilon }+\mu _{L}\rho )T\sqrt{%
T}e^{G_{\varepsilon }T}.  \notag
\end{eqnarray}%
We have to prove that the coefficients of $\rho $ and $E_{\varepsilon }$ are
positive, at least on a small interval. Since the approximating optimal
solution tends uniformly to the optimal solution in $(P),$ we infer that 
\begin{equation*}
G_{\varepsilon }\rightarrow G,\mbox{ }F_{i,\varepsilon }\rightarrow F_{i},%
\mbox{ as }\varepsilon \rightarrow 0\mbox{, }i=1,2,
\end{equation*}%
\begin{equation*}
F_{1}=l_{A}^{\ast }\left\Vert A^{\ast }\right\Vert _{\infty }+l_{I}^{\ast
}\left\Vert I^{\ast }\right\Vert +\mu _{L}\left\Vert L^{\ast }\right\Vert
_{\infty }+\mu _{L}y_{1},\mbox{ }F_{2}=\left\Vert A^{\ast }\right\Vert
_{\infty }+\left\Vert I^{\ast }\right\Vert _{\infty },
\end{equation*}%
\begin{equation*}
G=C\left( \left\Vert \beta _{A}\right\Vert _{\infty }+\left\Vert \beta
_{I}\right\Vert _{\infty }+\sigma +\mu _{A}+l_{A}^{\ast }+\mu
_{I}+l_{I}^{\ast }+\left\Vert \xi \right\Vert _{\infty }+l_{A}^{\ast
2}+l_{I}^{2}\right) ,
\end{equation*}%
with $C$ a constant. We note that the function $t\rightarrow 4\mu _{L}\alpha
_{0}e^{G_{\varepsilon }t}t\sqrt{t}$ is positive, strictly increasing for $%
t>0 $ and vanish at $0.$ We have 
\begin{equation}
4\mu _{L}\alpha _{0}e^{G_{\varepsilon }t}t\sqrt{t}\leq 4\mu _{L}\alpha
_{0}e^{Gt}t\sqrt{t}+O(\varepsilon )<1\mbox{ as }\varepsilon \rightarrow 0.
\label{41-0}
\end{equation}%
Then, there exists $T_{1}>0$ such that (\ref{41-0}) takes place on $%
(0,T_{1}).$

Next, we show that $\rho (1-4\mu _{L}\alpha _{0}e^{G_{\varepsilon }T}T\sqrt{T%
})-(F_{0}+4\alpha _{0}F_{1,\varepsilon }T\sqrt{T}e^{G_{\varepsilon }T})>0$
for $t\in (0,T_{2}).$ Indeed, this comes to 
\begin{eqnarray}
&&F_{0}+4\alpha _{0}(F_{1,\varepsilon }+\rho \mu _{L})t\sqrt{t}%
e^{G_{\varepsilon }t}  \label{41-1} \\
&\leq &F_{0}+4\alpha _{0}(F_{1,\varepsilon }+\rho \mu _{L})t\sqrt{t}%
e^{Gt}+O(\varepsilon )<\rho .  \notag
\end{eqnarray}%
The function $t\rightarrow F_{0}+4\alpha _{0}(F_{1,\varepsilon }+\rho \mu
_{L})t\sqrt{t}e^{Gt}$ is positive, strictly increasing and equal with $F_{0}$
at $0$ and we note that 
\begin{equation}
F_{0}=|L_{0}-y_{1}|<\rho .  \label{45}
\end{equation}%
This implies that there exists $T_{2}>0$ such that (\ref{41-1}) is satisfied
on $(0,T_{2}).$ It remains to check (\ref{45}).

We recall that we let $y_{1}\in (0,L_{0})$ and if we choose $\rho \in
(L_{0}-y_{1},\widehat{L}-y_{1})$ it follows that (\ref{45}), $y_{1}+\rho
\leq \widehat{L}$ and $y_{1}-\rho \leq \widehat{L}$ hold. The last one is
true because by the choice of $\rho $ we have $y_{1}-\rho \leq
2y_{1}-L_{0}\leq L_{0}<\widehat{L}.$ By (\ref{41})-(\ref{41-1}) we can write
that 
\begin{eqnarray}
\int_{0}^{T}\left\vert e_{\varepsilon }^{\prime }(t)\right\vert dt &\leq &%
\frac{\alpha _{0}F_{2,\varepsilon }(F_{1,\varepsilon }+\mu _{L}\rho )T\sqrt{T%
}e^{G_{\varepsilon }T}}{\rho (1-4\mu _{L}\alpha _{0}e^{G_{\varepsilon }T}T%
\sqrt{T})-(F_{0}+4\alpha _{0}F_{1,\varepsilon }T\sqrt{T}e^{G_{\varepsilon
}T})},  \notag \\
\mbox{for }T &\in &(0,T_{loc}),\mbox{ }T_{loc}:=\min \{T_{1},T_{2}\}.
\label{46}
\end{eqnarray}%
More exactly, $T_{1}$ and $T_{2}$ are the solutions to 
\begin{equation}
4\mu _{L}\alpha _{0}e^{Gt}t\sqrt{t}=1\mbox{ and }F_{0}+4\alpha
_{0}(F_{1,\varepsilon }+\rho \mu _{L})t\sqrt{t}e^{Gt}=\rho ,  \label{48}
\end{equation}%
respectively. Since the right-hand side in (\ref{46}) is bounded we conclude
with (\ref{29-2}), while (\ref{29-3}) is implied by (\ref{32}).\hfill $%
\square $

\medskip

\noindent \textbf{Remark 2.7. }Now, we recall a few definitions and results
necessary in the proof of the next theorem. We denote by $BV([0,T])$ the
space of functions $v:[0,T]\rightarrow \mathbb{R}$ with bounded variation,
that is 
\begin{equation*}
\left\Vert v\right\Vert _{BV([0,T])}=\sup \left\{
\sum\limits_{i=0}^{M-1}\left\vert v(t_{i+1})-v(t_{i})\right\vert ;\mbox{ }%
0=t_{0}<t_{1}<...<t_{M}=T\right\} <\infty ,
\end{equation*}%
and by $\mathcal{M}([0,T])$ the dual of the separable space $C([0,T]).$ The
space $\mathcal{M}([0,T])$ contains the bounded\ Radon measures defined on $%
[0,T].$ We also recall by the Lebesgue decomposition theorem (see e.g. \cite%
{Rockafeller-II-71}), that every $\mu \in (L^{\infty }(0,T))^{\ast }$ can be
uniquely written as 
\begin{equation}
\mu =\mu _{a}+\mu _{s},  \label{L1-2}
\end{equation}%
where $\mu _{a}\in L^{1}(0,T)$ and $\mu _{s}$ is a singular measure. This
means that for each $\varepsilon >0$ there exists a Lebesgue measurable set $%
\mathcal{S}\subset \lbrack 0,T]$ with meas$([0,T]\backslash \mathcal{S})\leq
\varepsilon $ and $\mu _{s}(\varphi )=0$ for all $\varphi \in L^{\infty }(%
\mathcal{S}).$ The support ($supp$) of $\mu _{s}$ is the set of all $t\in
\lbrack 0,T]$ for which $\mu _{s}(\varphi )\neq 0,$ for all $\varphi \in
L^{\infty }(t-\delta ,t+\delta ),$ and all $\delta $ positive.

Next, we recall that every $v\in BV([0,T])$ has a unique decomposition, 
\begin{equation}
v=v^{a}+v^{s},  \label{L1-3}
\end{equation}%
where $v^{a}\in AC[0,T]$ and $v^{s}\in BV([0,T]).$ Here, $AC[0,T]$ is the
space of absolutely continuous functions on $[0,T]$ and $v^{s}$ is a
singular part (for instance it can be a jump function with bounded variation
or a function with bounded variation with a.e. zero derivative).

\noindent We note that if $v\in BV([0,T]),$ then its distributional
derivative $\frac{dv}{dt}:=\mu $ belongs to $(L^{\infty }(0,T))^{\ast },$
and in virtue of the Lebesgue decomposition, it is represented by the sum of
the absolutely continuous part and the singular part%
\begin{equation}
\frac{dv}{dt}=\mu _{a}+\mu _{s}=\frac{dv^{a}}{dt}+\frac{dv^{s}}{dt}\in 
\mathcal{D}^{\prime }(0,T),  \label{L1-4}
\end{equation}%
where $\mathcal{D}^{\prime }(0,T)$ is the space of Schwartz distributions on 
$(0,T).$

\medskip

In the next theorem we shall pass to the limit in the approximating
optimality conditions. To this end, we introduce the system%
\begin{equation}
p^{\prime }-k_{0}^{\ast }p+k_{0}^{\ast }q=0,\mbox{ a.e. }t\in (0,T),
\label{17-1}
\end{equation}%
\begin{equation}
q^{\prime }-\beta _{A}S^{\ast }p+k_{3}^{\ast }q+\sigma d+l_{A}^{\ast }e+\mu
_{A}f=-\alpha _{0}A^{\ast },\mbox{ a.e. }t\in (0,T),  \label{18-1}
\end{equation}%
\begin{equation}
d^{\prime }-\beta _{I}S^{\ast }p+\beta _{I}S^{\ast }q-k_{2}^{\ast
}d+l_{I}^{\ast }e+\mu _{I}f=-\alpha _{0}I^{\ast },\mbox{ a.e. }t\in (0,T),
\label{19-1}
\end{equation}%
\begin{equation}
e^{\prime }-\mu _{L}e+\mu _{L}f+\nu =0,\mbox{ in }\mathcal{D}^{\prime }(0,T),
\label{20-1}
\end{equation}%
\begin{equation}
f^{\prime }+\xi p-\xi f=0,\mbox{ a.e. }t\in (0,T),  \label{22-1}
\end{equation}%
\begin{equation}
p(T)=0,\mbox{ }q(T)=0,\mbox{ }d(T)=0,\mbox{ }e(T)=0,\mbox{ }f(T)=0,
\label{23-1}
\end{equation}%
\begin{equation}
\nu =\alpha _{2}\lim_{\varepsilon \rightarrow 0}\left( \frac{1}{\varepsilon }%
(L_{\varepsilon }^{\ast }-\widehat{L})^{+}\right) \mbox{ weak* in }%
(L^{\infty }(0,T))^{\ast },  \label{23-2}
\end{equation}%
\begin{equation}
\nu ={\small \nu }_{a}+\nu _{s},\mbox{ }\nu _{a}\in L^{1}(0,T),\mbox{ }\nu
_{s}\in \mathcal{M}([0,T]),  \label{23-3}
\end{equation}%
where%
\begin{eqnarray}
k_{0}^{\ast } &=&\beta _{A}A^{\ast }+\beta _{I}I^{\ast },\mbox{ }k_{1}^{\ast
}=\sigma +\mu _{A}+l_{A}^{\ast },  \label{41-3} \\
k_{2}^{\ast } &=&\mu _{I}+l_{I}^{\ast },\mbox{ }k_{3}^{\ast }=\beta
_{A}S^{\ast }-k_{1}^{\ast }.  \notag
\end{eqnarray}%
We also define 
\begin{equation*}
\widehat{K}=\{y\in \mathbb{R};\mbox{ }-\infty <y\leq \widehat{L}\},\mbox{ }%
\widehat{\mathcal{K}}=\{v\in C([0,T]);\mbox{ }v(t)\leq \widehat{K}\mbox{ for
all }t\in \lbrack 0,T]\}.
\end{equation*}%
We denote by $N_{\widehat{K}}(z)$ the normal cone to $\widehat{K}$ at $z,$%
\begin{equation*}
N_{\widehat{K}}(z)=\left\{ 
\begin{array}{l}
0,\mbox{ \ \ \ \ if }z<\widehat{L} \\ 
\mathbb{R}^{+},\mbox{ \ \ if }z=\widehat{L}%
\end{array}%
\right. 
\end{equation*}%
and by%
\begin{equation}
\mathcal{N}_{\widehat{\mathcal{K}}}(\zeta )=\{\eta \in \mathcal{M}([0,T]);%
\mbox{ }\eta (\zeta -z)\geq 0,\mbox{ }\forall z\in \mbox{ }\widehat{\mathcal{%
K}}\}  \label{57}
\end{equation}%
the normal cone to $\widehat{\mathcal{K}}$ at $z,$ where $\eta (\zeta -z)$
is the value of the measure $\eta \in \mathcal{M}([0,T])$ at $(\zeta -z)\in
C([0,T]).$

\medskip

\noindent \textbf{Theorem 2.8. }\textit{Let }$(l_{A}^{\ast },l_{I}^{\ast })$%
\textit{\ be optimal in }$(P)$\textit{\ with the corresponding state }$%
X^{\ast }.$\textit{\ Then}, \textit{if }$T<T_{loc}$\textit{\ defined in} (%
\ref{46})-(\ref{48}), \textit{the optimality conditions for problem}\textbf{%
\ }$(P)$ \textit{read}:%
\begin{equation}
l_{A}^{\ast }=P_{[0,1]}\left( \frac{1}{\alpha _{1}}\int_{0}^{T}A^{\ast
}(q-e)dt\right) ,\mbox{ }l_{I}^{\ast }=P_{[0,1]}\left( \frac{1}{\alpha _{1}}%
\int_{0}^{T}I^{\ast }(d-e)dt\right)   \label{49}
\end{equation}%
\textit{where} $(p,q,d,e,f)$ \textit{is the solution of the dual system }(%
\ref{17-1})-(\ref{23-3}) \textit{with} 
\begin{equation}
(p,q,d,f)\in (W^{1,\infty }(0,T))^{4},\mbox{ }e\in BV([0,T],\mbox{ }%
e^{\prime }\in \mathcal{M}([0,T]).  \label{50-0}
\end{equation}%
\textit{Moreover, we have }$e^{\prime }=(e^{\prime })_{a}+(e^{\prime })_{s},$
$(e^{\prime })_{a}\in L^{1}(0,T),$ $(e^{\prime })_{s}\in \mathcal{M}([0,T]),$%
\begin{equation}
e_{a}^{\prime }(t)+\nu _{a}(t)=(\mu _{L}e-\mu _{L}f)(t),\mbox{ \textit{a.e.} 
}t\in (0,T),  \label{58}
\end{equation}%
\begin{equation}
(e_{s})^{\prime }+\nu _{s}=0,\mbox{ \ \ \ \ \ \textit{in} }\mathcal{D}%
^{\prime }(0,T),  \label{59}
\end{equation}%
\textit{where }%
\begin{equation}
\nu _{a}(t)\in N_{\widehat{K}}(L^{\ast }(t)),\mbox{ \textit{a.e.} }t\in
\lbrack 0,T],\mbox{ }  \label{60}
\end{equation}%
\begin{equation}
\nu _{s}(\varphi )\geq 0,\mbox{ }supp\mbox{ }\nu _{s}\subset \{t\in \lbrack
0,T];\mbox{ }L^{\ast }(t)=\widehat{L}\}.  \label{62}
\end{equation}

\medskip

\noindent \textbf{Proof. }We shall establish some estimates in order to pass
to the limit in the approximating optimality conditions determined in
Proposition 2.5. Let $T<T_{loc},$ let $(l_{A}^{\ast },l_{I}^{\ast })$\ be
optimal in $(P)$\ with the corresponding state $X^{\ast }$ and let us
consider $\{(l_{A,\varepsilon }^{\ast },l_{I,\varepsilon }^{\ast
}),X_{\varepsilon }^{\ast }\}$ optimal in $(P_{\varepsilon }).$ Recalling
Proposition 2.4 and the continuity property from Proposition 2.1, we have $%
l_{A,\varepsilon }^{\ast }\rightarrow l_{A}^{\ast },$ $l_{I,\varepsilon
}^{\ast }\rightarrow l_{I}^{\ast }$ and%
\begin{equation*}
X_{\varepsilon }^{\ast }\rightarrow X^{\ast }\mbox{ weak* in }W^{1,\infty
}(0,T)\mbox{ and uniformly on }[0,T].
\end{equation*}%
By (\ref{29-3}) each component of the solution to the dual system, but $%
e_{\varepsilon },$ follows to be bounded in $L^{\infty }(0,T).$ Moreover, by
(\ref{17})-(\ref{19}), (\ref{21}) it follows that $(p_{\varepsilon }^{\prime
})_{\varepsilon },(q_{\varepsilon }^{\prime })_{\varepsilon
},(d_{\varepsilon }^{\prime })_{\varepsilon },(f_{\varepsilon }^{\prime
})_{\varepsilon }$ are bounded in $L^{\infty }(0,T),$ so that, on a
subsequence, 
\begin{eqnarray}
p_{\varepsilon } &\rightarrow &p,\mbox{ }q_{\varepsilon }\rightarrow q,\mbox{
}d_{\varepsilon }\rightarrow d,\mbox{ }f_{\varepsilon }\rightarrow f\mbox{
weak* in }W^{1,\infty }(0,T)  \label{51} \\
&&\mbox{and uniformly in }[0,T],\mbox{ as }\varepsilon \rightarrow 0.  \notag
\end{eqnarray}%
The component $(e_{\varepsilon })_{\varepsilon }$ is bounded in $L^{\infty
}(0,T)$ by (\ref{38}) and its derivative is in $L^{1}(0,T),$ by (\ref{29-2}%
). These imply that $e_{\varepsilon }\in BV([0,T])$ and by Helly's theorem
(see e.g., \cite{vbp-2012}, p. 47) it follows that 
\begin{equation}
e_{\varepsilon }(t)\rightarrow e(t),\mbox{ for all }t\in \lbrack 0,T],\mbox{
as }\varepsilon \rightarrow 0.  \label{52}
\end{equation}%
Going back to (\ref{20}) we deduce that $\frac{1}{\varepsilon }%
(L_{\varepsilon }^{\ast }-\widehat{L})^{+}\in L^{1}(0,T).$ Now, we assert
that $(e_{\varepsilon }^{\prime })_{\varepsilon }$ and $\left( \frac{1}{%
\varepsilon }(L_{\varepsilon }^{\ast }-\widehat{L})^{+}\right) _{\varepsilon
}$ are weak* compact in $(L^{\infty }(0,T))^{\ast },$ the dual of $L^{\infty
}(0,T).$ This is pointed out in the proof of Corollary 2B in \cite%
{Rockafeller-II-71}, but this assertion does not follow directly from
Alaoglu theorem. An argument can be found in \cite{VB-DaPrato}, and we
resume it below.

Let us consider the linear operator $\Psi :C([0,T])\rightarrow L^{\infty
}(0,T),$ $\Psi z=\widetilde{\Psi },$ which maps a continuous function into
the corresponding class of equivalence $\widetilde{\Psi }$ (of all functions
a.e. equal). Its adjoint $\Psi ^{\ast }:(L^{\infty }(0,T))^{\ast
}\rightarrow \mathcal{M}([0,T]),$ is defined by $(\Psi ^{\ast }\mu )(z):=\mu
(\Psi z)$ for any $z\in C([0,T]).$ If $(\mu _{n})_{n}$ is bounded in $%
(L^{\infty }(0,T))^{\ast }$ and also in $\mathcal{M}([0,T]),$ then $(\Psi
^{\ast }\mu _{n})_{n}$ is bounded in $\mathcal{M}([0,T])$ and using the
Alaoglu theorem it follows that $(\Psi ^{\ast }\mu _{n})_{n}$ is weak*\
sequentially compact in $\mathcal{M}([0,T]).$ Therefore, it follows that $%
(\mu _{n})_{n}$ is weak*\ sequentially compact in $\mathcal{M}([0,T]).$
Passing to the limit in $\mu _{n}(\Psi z)=(\Psi ^{\ast }\mu _{n})(z)$ we get 
$\mu (\Psi z)=(\Psi ^{\ast }\mu )(z)$ for any $\widetilde{\Psi }\in
L^{\infty }(0,T)$ which is of the form $\Psi z$ with $z\in C([0,T]).$ Then,
due to the Hahn-Banach theorem, $\mu $ can be extended to all $L^{\infty
}(0,T)$ and so we conclude that $(\mu _{n})_{n}$ is weak*\ sequentially
compact in $(L^{\infty }(0,T))^{\ast }.$

Therefore, one can extract a subsequence such that 
\begin{equation}
e_{\varepsilon }^{\prime }\rightarrow q^{\prime },\mbox{ }\frac{\alpha _{2}}{%
\varepsilon }(L_{\varepsilon }^{\ast }-\widehat{L})^{+}\rightarrow \nu \mbox{
weak* in }(L^{\infty }(0,T))^{\ast }\subset \mathcal{M}([0,T]).  \label{53}
\end{equation}%
Thus, relying on (\ref{51})-(\ref{53}), we can pass to the limit in (\ref{17}%
)-(\ref{22}) and obtain (\ref{17-1})-(\ref{23-2}).

Now, we move to (\ref{23}), or more exactly in (\ref{28}) and (\ref{29-0})
and pass to the limit. The left-hand side of (\ref{28}) converges and since
the normal cone is maximal monotone, hence strongly-strongly closed it
follows that%
\begin{equation*}
-\int_{0}^{T}I^{\ast }(e-d)dt-\alpha _{1}l_{I}^{\ast }\in
N_{[0,1]}(l_{I}^{\ast }),
\end{equation*}%
which implies the second relation in (\ref{49}). Similarly, we proceed in (%
\ref{29-0}) and obtain 
\begin{equation}
-\int_{0}^{T}A^{\ast }(e-q)dt-\alpha _{1}l_{A}^{\ast }\in
N_{[0,1]}(l_{A}^{\ast }),  \notag
\end{equation}%
whence we get the first relation in (\ref{49}).

Finally, we detail equation (\ref{20-1}). Since $\frac{1}{\varepsilon }z^{+}$
is the subdifferential of the function $\frac{1}{2}\left( z^{+}\right) ^{2}$%
, we can write 
\begin{equation*}
\frac{\alpha _{2}}{\varepsilon }\int_{0}^{T}(L_{\varepsilon }^{\ast }-%
\widehat{L})^{+}((L_{\varepsilon }^{\ast }-\widehat{L})-(z(t)-\widehat{L}%
))dt\geq \frac{\alpha _{2}}{2\varepsilon }((L_{\varepsilon }^{\ast }-%
\widehat{L})^{+})^{2}-\frac{1}{2}((z(t)-\widehat{L})^{+})^{2}\geq 0\mbox{,}
\end{equation*}%
for all $z\in \widehat{\mathcal{K}}.$ At limit we obtain 
\begin{equation}
\nu (L^{\ast }-z)\geq 0,\mbox{ for all }z\in \mbox{ }\widehat{\mathcal{K}},
\label{56}
\end{equation}%
whence, $\nu \in \mathcal{N}_{\widehat{\mathcal{K}}}(L^{\ast }).$

Since $\nu \in (L^{\infty }(0,T))^{\ast }\subset \mathcal{M}([0,T]),$
recalling (\ref{L1-2}) and (\ref{L1-4}), we can represent $\nu =\nu _{a}+\nu
_{s},$ where $\nu _{a}$ is the absolutely continuous part (in the sense of
measure) and $\nu _{s}$ is the singular part of $\nu .$ Also, $e^{\prime
}\in \mathcal{M}([0,T])$ and $e^{\prime }=(e^{\prime })_{a}+(e^{\prime })_{s}
$ where $(e^{\prime })_{a}=e_{a}^{\prime }\in L^{1}(0,T)$. Then, (\ref{20-1}%
) can be rewritten as in (\ref{58})-(\ref{59}). Relation (\ref{56}) implies
that $\nu _{a}\in N_{\widehat{K}}(L^{\ast }),$%
\begin{equation}
\nu _{a}(t)\geq 0,\mbox{ }\nu _{a}(t)=0\mbox{ on }\{t\in \lbrack 0,T];\mbox{ 
}L^{\ast }(t)<\widehat{L}\},  \label{56-0}
\end{equation}%
\begin{equation}
\nu _{s}(\varphi )\geq 0,\mbox{ }\nu _{s}(\varphi )=0\mbox{ if }\varphi \in 
\overset{\circ }{\widehat{\mathcal{K}}},  \label{61}
\end{equation}%
where $\overset{\circ }{\widehat{\mathcal{K}}}=\{\varphi \in C([0,T]);$ $%
\varphi (t)<\widehat{L}$ for $t\in \lbrack 0,T]\}$ is the interior of $%
\widehat{\mathcal{K}},$ while $\nu _{s}$ has the support on the boundary of $%
\widehat{\mathcal{K}}.$ Recalling that 
\begin{equation*}
supp\mbox{ }\nu _{s}=\{\Sigma \subset \lbrack 0,T];\mbox{ }\nu _{s}\neq 0%
\mbox{ on }\Sigma \}
\end{equation*}%
it follows that $supp$ $\nu _{s}\subset \{t\in (0,T);$ $z(t)=\widehat{L}\}.$
Thus, we actually get (\ref{62}), as claimed. \hfill $\square $

\section{Problem $(P_{0})$}

\setcounter{equation}{0}

In this section we treat problem $(P_{0})$ associated to the cost functional
(\ref{10}).

\medskip

\noindent \textbf{Proposition 3.1. }\textit{Problem }$(P_{0})$\textit{\ has
at least one solution }$(\beta _{I}^{\ast },A_{0}^{\ast },I_{0}^{\ast })$.

\medskip

\noindent The proof is led on the basis of similar arguments as in
Proposition 2.2, using the result of existence and uniqueness of the
solution to the state system given in Proposition 2.1.

\medskip

The optimality conditions can be directly determined, after writing the
system in variations and the adjoint system. Let us define 
\begin{equation*}
K_{0}=\{(y,z)\in \mathbb{R}^{2};\mbox{ }y\geq 0,\mbox{ }z\geq 0,\mbox{ }%
y+z\leq N_{0}\},\mbox{ }K_{+}=\{y\in \mathbb{R};\mbox{ }y\geq 0\}.
\end{equation*}%
Let $\lambda >0$ and set the variations%
\begin{eqnarray}
\beta _{I}^{\lambda } &=&\beta _{I}^{\ast }+\lambda u,\mbox{ }u=\widetilde{%
\beta _{I}}-\beta _{I}^{\ast },\mbox{ where }\widetilde{\beta _{I}}(t)\in
K_{+},\mbox{ a.e. }t\geq 0,  \notag \\
A_{0}^{\lambda } &=&A_{0}^{\ast }+\lambda w,\mbox{ }w=\widetilde{A_{0}}%
-A_{0}^{\ast },\mbox{ }I_{0}^{\lambda }=I_{0}^{\ast }+\lambda v,\mbox{ }v=%
\widetilde{I_{0}}-I_{0}^{\ast },  \notag \\
\mbox{where }(\widetilde{A_{0}},\mbox{ }\widetilde{I_{0}}) &\in &K_{0}.
\label{201}
\end{eqnarray}%
Let us denote $x^{\lambda }:=\frac{X^{\beta _{I}^{\lambda },A_{0}^{\lambda
},I_{0}^{\lambda }}-X^{\ast }}{\lambda },$ where $X^{\beta _{I}^{\lambda
},A_{0}^{\lambda },I_{0}^{\lambda }}$ is the solution to (\ref{1})-(\ref{6})
corresponding to $(\beta _{I}^{\lambda },A_{0}^{\lambda },I_{0}^{\lambda }),$
satisfying (\ref{8}) and $X^{\ast }$ is the optimal state corresponding to $%
(\beta _{I}^{\ast },A_{0}^{\ast },I_{0}^{\ast }).$

We introduce the linearized system for problem $(P_{0})$%
\begin{eqnarray}
s^{\prime } &=&-k_{0}^{\ast }s-\beta _{I}^{\ast }S^{\ast }i-\beta
_{A}S^{\ast }a+\xi r-uS^{\ast }I^{\ast },  \notag \\
a^{\prime } &=&k_{0}^{\ast }s+\beta _{I}^{\ast }S^{\ast }i+k_{3}^{\ast
}a+uS^{\ast }I^{\ast },  \notag \\
i^{\prime } &=&\sigma a-k_{2}i,  \label{202} \\
l^{\prime } &=&l_{A}a+l_{I}i-\mu _{L}l,  \notag \\
r^{\prime } &=&\mu _{A}a+\mu _{I}i+\mu _{L}l-\xi r  \notag
\end{eqnarray}%
for $t\geq 0,$ with the initial condition 
\begin{equation}
s(0)=-w-v,\mbox{ }a(0)=w,\mbox{ }i(0)=v,\mbox{ }l(0)=0,\mbox{ }r(0)=0,
\label{203}
\end{equation}%
where 
\begin{equation}
k_{0}^{\ast }=\beta _{I}^{\ast }I^{\ast }+\beta _{A}A^{\ast },\mbox{ }%
k_{3}^{\ast }=\beta _{A}S^{\ast }-k_{1},\mbox{ }k_{1}=\sigma +\mu _{A}+l_{A},%
\mbox{ }k_{2}=\mu _{I}+l_{I}.  \label{204}
\end{equation}%
\noindent We introduce the backward dual system for the variables $%
(p,q,d,e,f)$ as 
\begin{equation}
p^{\prime }-k_{0}^{\ast }p+k_{0}^{\ast }q=0,  \label{205}
\end{equation}%
\begin{equation}
q^{\prime }-\beta _{A}S^{\ast }p+k_{3}^{\ast }q+\sigma d+l_{A}e+\mu _{A}f=0,
\label{206}
\end{equation}%
\begin{equation}
d^{\prime }-\beta _{I}^{\ast }S^{\ast }p+\beta _{I}^{\ast }S^{\ast
}q-k_{2}d+l_{I}e+\mu _{I}f=0,  \label{207}
\end{equation}%
\begin{equation}
e^{\prime }-\mu _{L}e+\mu _{L}f=0,  \label{208}
\end{equation}%
\begin{equation}
f^{\prime }+\xi p-\xi f=0,  \label{209}
\end{equation}%
for a.a. $t>0,$ with the final conditions 
\begin{eqnarray}
p(T) &=&0,\mbox{ }q(T)=0,\mbox{ }d(T)=0,\mbox{ }  \label{210} \\
e(T) &=&L^{\ast }(T)-L_{T},\mbox{ }f(T)=R^{\ast }(T)-R_{T}.  \notag
\end{eqnarray}%
The linear systems (\ref{202}) and (\ref{205})-(\ref{210}) have unique
global solutions in $(W^{1,\infty }(0,T))^{5}.$

Let $N_{K_{0}}(\zeta _{1},\zeta _{2})$ be the normal cone to the set $%
K_{0}\subset \mathbb{R}^{2}.$ We recall that $N_{K_{0}}(\zeta _{1},\zeta
_{2})=\partial I_{K_{0}}(\zeta _{1},\zeta _{2}),$ where $\partial
I_{K_{0}}(\zeta _{1},\zeta _{2})$ is the subdifferential of the indicator
set of $K_{0},$ that is 
\begin{equation*}
\partial I_{K_{0}}(\zeta _{1},\zeta _{2})=\{\eta =(\eta _{1},\eta _{2})\in 
\mathbb{R}^{2};\mbox{ }\eta _{1}(\zeta _{1}-z_{1})+\eta _{2}(\zeta
_{2}-\zeta _{2})\geq 0,\mbox{ for }(z_{1},z_{2})\in K_{0}\}.
\end{equation*}%
Let $N_{K_{+}}(z)$ be the normal cone to the set $K_{+},$ that is 
\begin{equation*}
N_{K_{+}}(z)=\left\{ 
\begin{array}{c}
\mathbb{R}^{-},\mbox{ if }z=0 \\ 
0,\mbox{ \ \ \ if }z>0.%
\end{array}%
\right.
\end{equation*}

\medskip

\noindent \textbf{Proposition 3.2. }\textit{Let }$(\beta _{I}^{\ast
},A_{0}^{\ast },I_{0}^{\ast })$\textit{\ be optimal in }$(P_{0})$ \textit{%
with the state }$X^{\ast }.$\textit{\ Then, }%
\begin{equation}
\beta _{I}^{\ast }(t)=P_{K_{+}}\left( \frac{1}{\alpha _{1}}%
(p(t)-q(t))S^{\ast }(t)I^{\ast }(t)\right) ,\mbox{ \textit{a.e.} }t\in (0,T),
\label{211}
\end{equation}%
\begin{equation}
(A_{0}^{\ast },I_{0}^{\ast })=(\Gamma +N_{K_{0}})^{-1}\left( \alpha
_{0}^{-1}(p(0)-q(0)+\alpha _{0}N_{0},p(0)-d(0)+\alpha _{0}N_{0})\right) ,
\label{212}
\end{equation}%
\textit{where }$\Gamma =\left( 
\begin{array}{cc}
2 & 1 \\ 
1 & 2%
\end{array}%
\right) $ \textit{and} $(p,q,d,e,_{\varepsilon })$\textit{\ is the solution
to the backward dual system} (\ref{205})-(\ref{210}).

\medskip

\noindent \textbf{Proof. }Let us multiply the equations for $s,a,i,l,r$ in (%
\ref{202}) by $p_{\varepsilon },q_{\varepsilon },d_{\varepsilon
},e_{\varepsilon },f_{\varepsilon },$ respectively and integrate over $%
(0,T). $ By integrating by parts we obtain%
\begin{eqnarray}
&&(L^{\ast }(T)-L_{T})l(T)+(R^{\ast }(T)-R_{T})r(T)  \label{213} \\
&=&\int_{0}^{T}S^{\ast }I^{\ast }(q-p)udt+p(0)(-w-v)+q(0)w+d(0)v.  \notag
\end{eqnarray}%
Since $(\beta _{I}^{\ast },A_{0}^{\ast },I_{0}^{\ast })$ is optimal in $%
(P_{0})$ we have $J(\beta _{I}^{\lambda },A_{0}^{\lambda },I_{0}^{\lambda
})\geq J(\beta _{I}^{\ast },A_{0}^{\ast },I_{0}^{\ast })$ and we deduce 
\begin{eqnarray}
&&(L^{\ast }(T)-L_{T})l(T)+(R^{\ast }(T)-R_{T})r(T)+\alpha
_{1}\int_{0}^{T}\beta _{I}^{\ast }udt  \label{214} \\
&&+\alpha _{0}\left( A_{0}^{\ast }w+I_{0}^{\ast }v+(N_{0}-A_{0}^{\ast
}-I_{0}^{\ast })(-w-v)\right) \geq 0.  \notag
\end{eqnarray}%
By comparison with (\ref{213}) and recalling the setting of $u,w,v,$ we have 
\begin{eqnarray}
&&\int_{0}^{T}[(p-q)S^{\ast }I^{\ast }-\alpha _{1}\beta _{I}^{\ast }](\beta
_{I}^{\ast }-\widetilde{\beta _{I}})dt  \label{215} \\
&&+[p(0)-q(0)-2\alpha _{0}A_{0}^{\ast }+\alpha _{0}N_{0}-\alpha
_{0}I_{0}^{\ast }](A_{0}^{\ast }-\widetilde{A_{0}})  \notag \\
&&+[p(0)-d(0)-2\alpha _{0}I_{0}^{\ast }+\alpha _{0}N_{0}-\alpha
_{0}A_{0}^{\ast }](I_{0}^{\ast }-\widetilde{I_{0}})\geq 0  \notag
\end{eqnarray}%
for all $\widetilde{\beta _{I}}(t)\in K_{+},$ a.e. $t\in (0,T)$ and $(%
\widetilde{A_{0}},\widetilde{I_{0}})\in K_{0}.$ In particular, by setting $(%
\widetilde{A_{0}},\widetilde{I_{0}})=(A_{0}^{\ast },I_{0}^{\ast }),$ (\ref%
{215}) yields 
\begin{equation}
(p(t)-q(t))S^{\ast }(t)I^{\ast }(t)-\alpha _{1}\beta _{I}^{\ast }(t)\in
N_{K_{+}}(\beta _{I}^{\ast }(t)),\mbox{ a.e. }t\in (0,T),  \label{28-00}
\end{equation}%
which implies (\ref{211}).

Then, let us set $\widetilde{\beta _{I}}=\beta _{I}^{\ast }.$ By (\ref{215})
we obtain that 
\begin{eqnarray*}
&&(p(0)-q(0)-2\alpha _{0}A_{0}^{\ast }+\alpha _{0}N_{0}-\alpha
_{0}I_{0}^{\ast },p(0)-d(0)-2\alpha _{0}I_{0}^{\ast }+\alpha
_{0}N_{0}-\alpha _{0}A_{0}^{\ast }) \\
&\in &N_{K_{0}}(A_{0}^{\ast },I_{0}^{\ast }).
\end{eqnarray*}%
This can be still written 
\begin{equation}
\alpha _{0}^{-1}(p(0)-q(0)+\alpha _{0}N_{0},p(0)-d(0)+\alpha
_{0}N_{0})-\Gamma (A_{0}^{\ast },I_{0}^{\ast })\in N_{K_{0}}(A_{0}^{\ast
},I_{0}^{\ast }),  \label{29-00}
\end{equation}%
where $\Gamma =\left( 
\begin{array}{cc}
2 & 1 \\ 
1 & 2%
\end{array}%
\right) .$ This implies (\ref{212}), because $N_{K_{0}}$ is maximal monotone
in $\mathbb{R}^{2}$, and $\Gamma :\mathbb{R}^{2}\rightarrow \mathbb{R}^{2}$
is positive, so $(\Gamma +N_{K_{0}})^{-1}$ is Lipschitz.\hfill $\square $

\medskip

\noindent \textbf{Remark 3.3. }By (\ref{211}) and the fact that the optimal
state is nonnegative, it is clear that\textbf{\ }%
\begin{equation}
\beta _{I}^{\ast }(t)=\left\{ 
\begin{tabular}{ll}
$0$ & $\mbox{on }\{t\in \lbrack 0,T];\mbox{ }p(t)\leq q(t)\}$ \\ 
$\frac{1}{\alpha _{1}}(p(t)-q(t))S^{\ast }(t)I^{\ast }(t)$ & on $\{t\in
\lbrack 0,T];\mbox{ }p(t)>q(t)\}.$%
\end{tabular}%
\right.   \label{63}
\end{equation}

\section{System stability and determination of the reproduction rate}

\setcounter{equation}{0}

In this section we investigate the system stability, which will help to
derive an expression for the reproduction rate. Once identified $\beta
_{I}(t)$, one can consider its average over $(0,T)$ and use it for the
stability analysis. An average can be set for $\beta _{A},$ too. The use of
the average can be more accurate if $(0,T)$ is short or if the epidemic has
reached a plateau, where the rates do not have large variations.

We discuss the system stability for $\xi =0,$ because the situation with $%
\xi $ constant nonzero$,$ meaning that immunity is lost immediately after
recovery is not realistic. We recall that $N=1.$

\medskip

\noindent \textbf{Theorem 4.1.}\textit{\ Let }$\xi =0$ \textit{and assume
that} 
\begin{equation}
\overline{S_{\infty }}:=\frac{(\sigma +\mu _{A}+l_{A})(\mu _{I}+l_{I})}{%
\beta }\leq N,\mbox{ }\beta :=k_{2}\beta _{A}+\sigma \beta _{I}.
\label{77-0}
\end{equation}%
\textit{The system }$(A,I,L)$ \textit{with a positive susceptible population}
$S_{\infty }$ is \textit{asymptotically stable if and only if } 
\begin{equation}
S_{\infty }<\overline{S_{\infty }}.  \label{77}
\end{equation}%
\textit{Moreover, all solutions (starting from any nonnegative initial
condition) tend to a stationary state, that is }%
\begin{equation}
\lim_{t\rightarrow \infty }(S(t),A(t),I(t),L(t),R(t))=(\widetilde{S_{\infty }%
},0,0,0,N-\widetilde{S_{\infty }})  \label{77-1}
\end{equation}%
\textit{exists and} $\widetilde{S_{\infty }}<\overline{S_{\infty }}.$

\medskip

\noindent \textbf{Proof.} If $\xi =0,$ the stationary solutions are found as 
$(S_{\infty },0,0,0,R_{\infty }),$ and we choose $S_{\infty },$ $R_{\infty
}\in \lbrack 0,N].$ We consider the linearized system, extract the system
for the infected compartments $A,$ $I$, $L,$ and define its matrix 
\begin{equation*}
\mathcal{A}_{1,lin}=\left( 
\begin{array}{ccc}
k_{3,\infty } & \beta _{I}S_{\infty } & 0 \\ 
\sigma & -k_{2} & 0 \\ 
l_{A} & l_{I} & -\mu _{L}%
\end{array}%
\right) ,
\end{equation*}
with $k_{1},$ $k_{2}$ given in (\ref{11}) and $k_{3,\infty }=\beta
_{A}S_{\infty }-k_{1}.$ The characteristic equation has a negative solution $%
\lambda =-\mu _{L}$ and two solutions to the equation 
\begin{equation}
P(\lambda )=\lambda ^{2}+\lambda (k_{1}+k_{2}-\beta _{A}S_{\infty
})+k_{1}k_{2}-\beta S_{\infty }=0.\mbox{ }  \label{84}
\end{equation}%
We prove that the solutions to (\ref{84}) have the real part negative,
meaning that the polynomial $P(\lambda )$ is Hurwitz.

The condition to be Hurwitz is that $k_{1}+k_{2}-\beta _{A}S_{\infty }>0$
and $k_{1}k_{2}-\beta S_{\infty }>0,$ that is $S_{\infty }\in \left[ 0,\min
\left\{ \frac{k_{1}k_{2}}{\beta },\frac{k_{1}+k_{2}}{\beta _{A}}\right\}
\right) =\left[ 0,\frac{k_{1}k_{2}}{\beta }\right) \subset \lbrack 0,N)$ by (%
\ref{77-0}) and this leads to (\ref{77}).

Let us prove (\ref{77-1}). By setting $\omega =\left( 
\begin{array}{ccc}
1 & 0 & 0%
\end{array}%
\right) ^{\mbox{T}}$ and denoting $Z:=(A,I,L)^{\mbox{T}}$ ($^{\mbox{T}}$ is
the transposed) we write by (\ref{1})-(\ref{5}) the equations%
\begin{eqnarray}
Z^{\prime }(t) &=&\mathcal{A}Z(t)+\omega (\beta _{A}A+\beta _{I}I)S(t)
\label{81} \\
&=&\left( 
\begin{array}{ccc}
-k_{1} & 0 & 0 \\ 
\sigma & -k_{2} & 0 \\ 
l_{A} & l_{I} & -\mu _{L}%
\end{array}%
\right) Z(t)+\left( 
\begin{array}{c}
1 \\ 
0 \\ 
0%
\end{array}%
\right) (\beta _{A}A+\beta _{I}I)S(t),  \notag
\end{eqnarray}%
\begin{equation}
S^{\prime }(t)=-S(t)(\beta _{A}A+\beta _{I}I).  \label{82}
\end{equation}%
We must show that there exists 
\begin{equation}
\lim_{t\rightarrow \infty }S(t)=\widetilde{S_{\infty }},\mbox{ }%
\lim_{t\rightarrow \infty }Z(t)=(0,0,0).  \label{83}
\end{equation}%
Since by (\ref{82}) we see that $S^{\prime }(t)\leq 0$ it follows that $%
S^{\prime }$ is monotonically decreasing and so it tends to a limit $%
\widetilde{S_{\infty }}\geq 0.$ Then, by (\ref{81})-(\ref{82}) we have 
\begin{equation*}
Z^{\prime }(t)=\mathcal{A}Z(t)-S^{\prime }(t)\omega
\end{equation*}%
and deduce by the formula of variation of constants and integration by
parts, that%
\begin{eqnarray}
Z(t) &=&e^{\mathcal{A}t}Z(0)-\int_{0}^{t}e^{\mathcal{A}(t-s)}S^{\prime
}(s)\omega ds  \label{86} \\
&=&e^{\mathcal{A}t}Z(0)-S(t)\omega +e^{\mathcal{A}t}S(0)\omega -\int_{0}^{t}%
\mathcal{A}e^{\mathcal{A}(t-s)}S(s)\omega ds.  \notag
\end{eqnarray}%
We calculate the last term%
\begin{eqnarray*}
&&E=\int_{0}^{t}\mathcal{A}e^{\mathcal{A}(t-s)}S(s)\omega ds=E_{1}+E_{2} \\
E_{1} &=&\int_{0}^{t}\mathcal{A}e^{\mathcal{A}(t-s)}\widetilde{S_{\infty }}%
\omega ds,\mbox{ }E_{2}=\int_{0}^{t}\mathcal{A}e^{\mathcal{A}(t-s)}(S(s)-%
\widetilde{S_{\infty }})\omega ds.
\end{eqnarray*}%
and get 
\begin{equation*}
E_{1}=-\int_{0}^{t}\mathcal{(}e^{\mathcal{A}(t-s)})^{\prime }\widetilde{%
S_{\infty }}\omega ds=-\mathcal{(}1-e^{-\mathcal{A}t})\widetilde{S_{\infty }}%
\omega =-\widetilde{S_{\infty }}\omega +e^{\mathcal{A}t}\widetilde{S_{\infty
}}\omega .
\end{equation*}%
We note that $\mathcal{A}$ is Hurwitz and recall that $\left\Vert e^{%
\mathcal{A}t}\right\Vert \leq Me^{-\gamma t},$ where $\gamma <\min_{j}\{-$Re(%
$\lambda _{j});$ $\lambda _{j}$ are the eigenvalues of $\mathcal{A}\}$ and $%
\left\Vert \cdot \right\Vert $ is the norm in $\mathbb{R}^{3}$.

Now, for any $M,$ such that $t>M$ we write%
\begin{eqnarray*}
E_{2} &=&\mathcal{A}\int_{0}^{t}e^{\mathcal{A}(t-s)}(S(s)-\widetilde{%
S_{\infty }})\omega ds \\
&=&\mathcal{A}\int_{0}^{M}e^{\mathcal{A}(t-s)}(S(s)-\widetilde{S_{\infty }}%
)\omega ds+\mathcal{A}\int_{M}^{t}e^{\mathcal{A}(t-s)}(S(s)-\widetilde{%
S_{\infty }})\omega ds \\
&\leq &C\sup_{0\leq s\leq M}\left\vert (S(s)-\widetilde{S_{\infty }})\omega
\right\vert \left\Vert \mathcal{A}\right\Vert \int_{0}^{M}e^{-\mathcal{%
\gamma (}t-s)}ds \\
&&+C\sup_{M\leq s<\infty }\left\vert (S(s)-\widetilde{S_{\infty }})\omega
\right\vert \left\Vert \mathcal{A}\right\Vert \int_{M}^{t}e^{-\mathcal{%
\gamma (}t-s)}ds \\
&\leq &\frac{C}{\gamma }\sup_{0\leq s\leq M}\left\vert (S(s)-\widetilde{%
S_{\infty }})\omega \right\vert \left\Vert \mathcal{A}\right\Vert e^{-\gamma
t}[e^{\gamma M}-1] \\
&&+\frac{C}{\gamma }\sup_{M\leq s<\infty }\left\vert (S(s)-\widetilde{%
S_{\infty }})\omega \right\vert \left\Vert \mathcal{A}\right\Vert
[1-e^{-\gamma (t-M)}],
\end{eqnarray*}%
where $\left\Vert \mathcal{A}\right\Vert $ is the norm of the matrix $%
\mathcal{A}.$ The first term of the last sum tends to zero as $t\rightarrow
\infty .$ For the second we give the following argument. Let $\varepsilon >0$
and fix $M$ such that $\sup_{M\leq s<\infty }\left\vert (S(s)-\widetilde{%
S_{\infty }})\omega \right\vert <\varepsilon .$ Therefore, 
\begin{equation*}
\limsup_{t\rightarrow \infty }\sup_{M\leq s<\infty }\left\vert (S(s)-%
\widetilde{S_{\infty }})\omega \right\vert [1-e^{-\gamma (t-M)}]<\varepsilon
\end{equation*}%
and since $\varepsilon $ is arbitrary, it means that this term tends to
zero, too. Thus, $\lim_{t\rightarrow \infty }E_{2}(t)=0.$ We note that 
\begin{equation*}
\lim_{t\rightarrow \infty }e^{\mathcal{A}t}Z(0)=0,\mbox{ }\lim_{t\rightarrow
\infty }e^{\mathcal{A}t}S(0)\omega =0,\mbox{ }\lim_{t\rightarrow \infty }e^{%
\mathcal{A}t}\widetilde{S_{\infty }}\omega =0,
\end{equation*}%
since $\mathcal{A}$ is Hurwitz. Then, taking into account the first relation
in (\ref{83}) and letting $t\rightarrow \infty $ in (\ref{86}) we get 
\begin{equation*}
\lim_{t\rightarrow \infty }Z(t)=-\widetilde{S_{\infty }}\omega +\widetilde{%
S_{\infty }}\omega =0.
\end{equation*}%
Thus, (\ref{77-1}) follows and (\ref{77}) is a sufficient condition to have
the solution $(A,I,L)=(0,0,0)$ stable.

It remains to show that $\widetilde{S_{\infty }}<\overline{S_{\infty }}.$
Otherwise, if $S_{\infty }\geq \overline{S_{\infty }},$ we see\ that the
polynomial is no longer Hurwitz and so system $(A,I,L)$ is not
asymptotically stable. Moreover, $\lim_{t\rightarrow \infty
}R(t)=N-(S(t)+A(t)+I(t)+L(t))=N-S_{\infty }:=R_{\infty }$ and we note that $%
R_{\infty }=N-S_{\infty }>N-\overline{S_{\infty }}=N-\frac{k_{1}k_{2}}{\beta 
}\geq 0,$ according to (\ref{77-0}).

We prove that condition (\ref{77}) is necessary, too. It means that if $%
S_{\infty }$ is a steady state and $\lim_{t\rightarrow \infty
}(A(t),I(t),L(t))=(0,0,0)$, it follows that $S_{\infty }$ should satisfy (%
\ref{77}). Let us assume the opposite, that is $S_{\infty }\geq \overline{%
S_{\infty }}.$ It follows that the above corresponding polynomial $P(\lambda
)$ is not Hurwitz and this implies that the system $(A,I,L)$ is not
asymptotically stable, meaning that $(A(t),I(t),L(t))$ does no longer tend
to $0$ as $t\rightarrow \infty $.\hfill $\square $

\bigskip

The fact that the value $\frac{\beta }{k_{1}k_{2}}$ appears as a critical
value for the system, gives a justification to define the reproduction rate
as 
\begin{equation}
R^{0}=\frac{k_{2}\beta _{A}+\sigma \beta _{I}}{(\sigma +\mu _{A}+l_{A})(\mu
_{I}+l_{I})}.  \label{70-1}
\end{equation}%
Its epidemiologic interpretation will be given further.

\section{Conclusions}

\setcounter{equation}{0}

We solved an inverse and a control problem related to an epidemic model for
SARS-CoV-2, with five compartments: susceptible $S,$ undetected infected
asymptomatic $A,$ undetected infected symptomatic $I$, detected by testing
and isolated $L,$ and recovered $R$. By an optimal control technique we
identified the rate of infection $\beta _{I}^{\ast }$ of the susceptible
individuals by the infective $I$ class and the number of undetected
individuals in the classes $A$ and $I$ at a time set $t=0.$ Their estimation
was relied on the observation of the number of the isolated and recovered
people at time $t=0$ and at another later time $T.$ It turned out that the $%
\beta _{I}^{\ast }$ is the projection of a point depending on the optimal
states and the solution to the dual system, on the set of positive real
numbers and the point $(A_{0}^{\ast },I_{0}^{\ast })$ is uniquely determined
by (\ref{212}). Knowing the evolution of $\beta _{I}^{\ast }$ during $(0,T)$
and the initial values $A_{0}^{\ast }$ and $I_{0}^{\ast }$ one can estimate
the transient evolution of the compartments after the time $T.$ Then, the
control of the infected classes $A$ and $I$ was done, by means of the
coefficients $l_{A}$ and $l_{I}$ related to the testing action. The
controllers were provided by expressions depending on the optimal states for
this problem and the solution to the singular dual system.

The investigation of the system asymptotic stability enhanced the
determination of the reproduction rate $R^{0},$ defined in (\ref{70-1}),
under the assumption of constant coefficients in the state system. Theorem
5.1 characterizes the system behavior under the assumption that the disease
induces life immunity. It indicates an asymptotic extinction of the disease,
following by the globally asymptotic stability of the solution to a steady
state $(\widetilde{S_{\infty }},0,0,0,N-S_{\infty }),$ where $\widetilde{%
S_{\infty }}$ does not exceed the value $\overline{S_{\infty }}$ given by (%
\ref{77}). This is interpreted as the number of individuals that have been
never infected (see \cite{Gio}). Moreover, if $R^{0}\widetilde{S_{\infty }}%
<1,$ the epidemic extinguishes, while the case $R^{0}\widetilde{S_{\infty }}%
>1$ shows a massive outbreak, $R^{0}\widetilde{S_{\infty }}=1$ being a
bifurcation point.

Finally, we underline that, in order to avoid much more calculations in a
model with many equations, we used a restraint model with less compartments,
including however the most relevant ones. More accurate values for the
desired parameters to be identified can be obtained developing similar
arguments for a more elaborated model with many compartments supposed to be
measurable, such that the information provided by their observation could be
included in the minimization problem formulation. Also, other parameters, as
for example $\beta _{A}$ can be identified and numerical simulations will be
provided in a forthcoming paper.

\end{document}